\documentclass{amsart}
\usepackage{amscd}
\newtheorem{thm}{Theorem}[section]
\newtheorem{lemma}[thm]{Lemma}
\newtheorem{prop}[thm]{Proposition}
\newtheorem{cor}[thm]{Corollary}

\newtheorem{conj}[thm]{Conjecture}
\newtheorem{remark}[thm]{Remark}
\theoremstyle{definition}
\newtheorem{defn}[thm]{Definition}
\newtheorem{exam}[thm]{Example}
\newtheorem{nota}[thm]{}
\numberwithin{equation}{section}

\newcommand{\Ass}{\mathrm{Ass}}
\newcommand{\Assh}{\mathrm{Assh}}
\newcommand{\chara}{\mathrm{char}\,}
\newcommand{\codim}{\mathrm{codim}\,}
\newcommand{\depth}{\mathrm{depth}}

\newcommand{\emb}{\mathrm{emb}\,}

\newcommand{\frm}{{\mathfrak m}}

\newcommand{\Hom}{\mathrm{Hom}}
\newcommand{\Ker}{\mathrm{Ker}\,}
\newcommand{\Image}{\mathrm{Im}\,}
\newcommand{\indeg}{\mathrm{indeg}\,}
\newcommand{\link}{\mathrm{link}}
\newcommand{\Mat}{\mathrm{Mat}\,}

\newcommand{\reg}{\mathrm{reg}\,}
\newcommand{\rank}{\mathrm{rank}\,}

\newcommand{\Tor}{\mathrm{Tor}}
\newcommand{\A}{\ensuremath{\mathbb A}}
\newcommand{\F}{\ensuremath{\mathbb F}}
\newcommand{\G}{\ensuremath{\mathbb G}}
\newcommand{\Z}{\ensuremath{\mathbb Z}}
\newcommand{\N}{\ensuremath{\mathbb N}}

\newcommand{\bbR}{\ensuremath{\mathbb R}}
\newcommand{\bbZ}{\ensuremath{\mathbb Z}}
\newcommand{\qle}{\;\le\;}

\newcommand{\covDelta}{\widetilde{\Delta}}
\newcommand{\ecm}{e_{\mathrm{CM}}}
\newcommand{\mucm}{\mu_{\mathrm{CM}}}
\newcommand{\bbinom}[2]{%
\genfrac{(}{)}{0pt}{}{#1}{#2}}

\newcommand{\sbinom}[2]{%
\text{\small $\displaystyle{\left(\!\!\!\begin{array}{c}#1 \\ #2 \end{array}\!\!\!
\right)}$}}

\begin{document}
\title{Buchsbaum Stanley--Reisner rings with minimal multiplicity}
\author{Naoki Terai}
\address{Department of Mathematics, Faculty of Culture and Education, 
         Saga University, Saga 840--8502, Japan} 
\email{terai@cc.saga-u.ac.jp}
\author{Ken-ichi Yoshida}
\address{Graduate School of Mathematics, Nagoya University,
         Nagoya  464--8602, Japan}
\email{yoshida@math.nagoya-u.ac.jp}
\date{Dec,25, 2003} 
\begin{abstract}
In this paper, we study non-Cohen--Macaulay Buchsbaum 
Stanley--Reisner rings with linear free resolution. 
In particular, for given integers $c$, $d$, $q$ with 
$c \ge 1$, $2 \le q \le d$, 
we give an upper bound $h_{c,d,q}$ on the dimension of 
the unique non-vanishing homology 
$\widetilde{H}_{q-2}(\Delta;k)$ of a $d$-dimensional Buchsbaum 
ring $k[\Delta]$ with $q$-linear resolution and codimension $c$. 
Also, 
we discuss about existence for such Buchsbaum rings 
with $\dim_k \widetilde{H}_{q-2}(\Delta;k) = h$ 
for any $h$ with $0 \le h \le h_{c,d,q}$, and prove an 
existence theorem  in the case of $q=d=3$
using the notion of Cohen--Macaulay linear cover. 
On the other hand, we introduce the notion of Buchsbaum Stanley--Reisner 
rings with minimal multiplicity of type $q$, 
which extends the notion of Buchsbaum rings with 
minimal multiplicity defined by Goto. 
As an application, we give many examples of Buchsbaum 
Stanley--Reisner rings with $q$-linear resolution. 
\end{abstract}
\maketitle
\section*{Introduction}
\par 
For any simplicial complex $\Delta$ on $V = \{x_1,\ldots,x_n\}$, 
the homogeneous reduced $k$-algebra 
$k[\Delta] = k[X_1,\ldots,X_n]/I_{\Delta}$, where $I_{\Delta}$ 
is the ideal generated by all square-free monomials 
$X_{i_1}\cdots X_{i_p}$ such that 
$\{x_{i_1},\ldots,x_{i_p}\} \notin \Delta$, 
is called the \textit{Stanley--Reisner ring} of $\Delta$. 
In recent years, the class of Stanley--Reisner rings is one 
of important classes in the theory of commutative algebra. 

\par \vspace{1mm}
In \cite{EiGo}, Eisenbud and Goto investigated rings with linear 
resolution and showed the significance of this property. 
Then from many viewpoints, they have widely studied. 
Let us pick up some important results in 
the class of Stanley--Reisner rings.    
Fr\"oberg \cite{Fr1, Fr2} classified all $\Delta$ for which $k[\Delta]$ 
has $2$-linear resolution. 
Hibi \cite{Hi2} gave a necessary and sufficient condition 
for a Buchsbaum Stanley--Reisner ring 
to have linear resolution in terms of the reduced homology 
of the simplicial complex and the $a$-invariants of its links. 
\par 
Also, there is a well-known criterion for a Cohen--Macaulay 
(Stanley--Reisner) ring to have linear resolution 
in terms of its $h$-vector or its multiplicity 
with given initial degree and codimension (see e.g. \cite{EiGo}). 
However, as for Buchsbaum case, it seems that there is no such a criterion. 
Hence, in this paper, we investigate 
the structure of Buchsbaum Stanley--Reisner rings with 
linear resolution in connection with their multiplicities.   

\par \vspace{1mm}
The purpose of this paper is divided into the following three pieces$:$
\begin{description}
 \item[(I)]  To give fundamental properties of Buchsbaum Stanley--Reisner 
rings with linear resolution. 
 \item[(II)] To introduce the notion of Buchsbaum Stanley--Reisner rings 
with minimal multiplicity of type $q$ for any integer $q \ge 2$. 
 \item[(III)] To construct $3$-dimensional Buchsbaum Stanley--Reisner rings 
having $3$-linear resolution with given parameters (codimension, 
dimension of the first reduced homology).    
\end{description}

\par \vspace{1mm}
Let us explain the organization of this paper. 
\par 
After recalling the notation and the terminology, 
in Section 2, we give fundamental properties of Buchsbaum Stanley--Reisner 
rings with linear resolution. 
Now let $A=k[\Delta]$ be a $d$-dimensional Buchsbaum Stanley--Reisner ring 
with $q$-linear resolution, and put $\codim A = c (\ge 1)$. 
Then $q \le d+1$. 
If $q=d+1$, then $\Delta$ is a $(d-1)$-skeleton of $2^{V}$ 
(Proposition \ref{d+1-linear}). 
So we may assume that $2 \le q \le d$.     
Then $H_{\frm}^{i}(A) =0$ for all $i\ne q-1$, $d$;  
$H_{\frm}^{q-1}(A) \cong \widetilde{H}_{q-2}(\Delta;k)$.  
Thus it seems that 
$h := \dim_k H_{\frm}^{q-1}(A)$ is an important 
invariant of $A$. 
From this point of view, we determined the $h$-vector of $A$ and 
proved an inequality$:$
\[
  0 \le h \le h_{c,d,q}:= 
\frac{(c+q-2)\cdots (c+1)c}{d(d-1)\cdots (d-q+2)}; 
\]
see Theorem \ref{Key}. Also, we pose the following conjecture$:$

\par \vspace{1mm} \par \noindent 
\textbf{Conjecture \ref{Conj-1}.}
 Let $d$, $c$, $q$, $h$ be integers with $c \ge 1$, $h \ge 0$, 
and $2 \le q \le d$. 
Then the following conditions  are equivalent$:$
\begin{enumerate}
\item There exists a Buchsbaum Stanley--Reisner ring $A=k[\Delta]$ 
with $q$-linear resolution
such that $\dim A =d$, $\codim A = c$ and 
$\dim_k H_{\frm}^{q-1}(A)=h$.
\item The above inequality $0 \le h \le h_{c,d,q}$ holds. 
\end{enumerate}

\par \vspace{1mm}
The coarse version, which we solve in Section 4, of this conjecture 
is given by Hibi in \cite{Hi2}. 

\par \vspace{1mm} \par \noindent 
{\bf Hibi's problem.}
Construct a Buchsbaum complex of dimension $d-1$ with $q$-linear 
resolution for any given integers $q,\,d$ with $2 \le q \le d$. 

\par \vspace{1mm}
In Section 3, we study the Alexander dual of Buchsbaum complexes  
with linear resolution. 
In fact, we prove that $k[\Delta]$ is Buchsbaum with $q$-linear resolution 
if and only if $k[\Delta^{*}]$ ($\Delta^{*}$ denotes 
the Alexander dual of $\Delta$)
has almost $c$-linear resolution 
and $\beta_{qj}(k[\Delta^{*}]) = 0$ for all 
$j \ne c+q-1$, $c+d$; see Theorem \ref{LinAlex}. 

\par \vspace{1mm}
In Section 4, 
we introduce the notion of 
minimal multiplicity of type $q$ 
for Buchsbaum Stanley--Reisner rings, and 
investigate its property. 
\par
In \cite{Go2}, Goto defined Buchsbaum local rings 
with minimal multiplicity, 
and  he proved that they have $2$-linear resolutions; 
see \cite{Go1, Go2}. 
We generalize this notion in the class of Stanley--Reisner rings
(Proposition \ref{GotoMinimal}).  
Namely, we prove the following theorem, which is a main result 
in this paper. 

\par \vspace{1mm} \par \noindent 
\textbf{Theorem \ref{Main}.}
Let $A = k[\Delta]$ be a Buchsbaum Stanley--Reisner ring 
with $\codim A = c$ and $\indeg A = q$. 
Then 
\begin{enumerate}
 \item The following inequality holds$:$
\[
   e(A) \ge \frac{c+d}{d} \bbinom{c+q-2}{q-2}. 
\]
 \item If the equality holds in $(1)$,  
 then it has  $q$-linear resolution. 
\end{enumerate}

\par \vspace{1mm}
We say that $A$ has \textit{minimal multiplicity of type $q$} if 
one of the following equivalent conditions
 (see Theorem \ref{Char} for more details)$:$  
\begin{enumerate}
 \item The equality holds in Theorem \ref{Main}(1).  
 \item $A$ has  $q$-linear resolution and 
$\dim_k H_{\frm}^{q-1}(A) = h_{c,d,q}$. 
 \item $a(A) = q-d-2$. 
 \item $k[\Delta^{*}]$ is Cohen--Macaulay with pure and almost linear 
 resolution with $a$-invariant $0$.  
\end{enumerate}

\par \vspace{1mm}
In the proof of Theorem \ref{Main}, 
\textit{Hochster's formula} and 
\textit{Hoa--Miyazaki theorem} play important 
roles. 
Also, using the same idea of the proof of 
Theorem \ref{Main}, we can 
prove a generalization of \textit{Hibi's criterion},  
which asserts that  
$A = k[\Delta]$ has  $q$-linear resolution if and only if 
$\widetilde{H}_{q-1}(\Delta;k) =0$ and 
$a(k[\link_{\Delta} \{x_i\}) \le q-d$ for all $i=1,\ldots,n$; 
see Theorem \ref{AinLink}.  
\par 
On the other hand, as an application of Theorems \ref{LinAlex} 
and \ref{Char}, we prove that the Alexander dual of the boundary 
complex of a cyclic polytope is Buchsbaum with minimal 
multiplicity in our sense. 
In particular, this example gives an affirmative answer to 
the above Hibi's problem. 
See Section 4 for more examples. 

\par \vspace{1mm}
In Section 5, we introduce the notion of Cohen--Macaulay 
cover, and prove that Conjecture \ref{Conj-1} is true for $q=d=3$. 
Let $\Delta$ be a pure simplicial complex on vertex set $V$ with 
$\indeg k[\Delta] = \dim k[\Delta] = d$. 
Then $\widetilde{\Delta}$ is said to be a 
\textit{Cohen--Macaulay cover} of $\Delta$ over a field $k$ 
if it is a $(d-1)$-dimensional simplicial complex on the same vertex 
set $V$ as containing $\Delta$ and 
$k[\widetilde{\Delta}]$ is Cohen--Macaulay with $d$-linear 
resolution. 
The notion of Cohen--Macaulay cover is very useful to attack 
the above conjecture. 
\par
Now consider the case of $q=d$ in the above Conjecture \ref{Conj-1}. 
Let $c$, $d$ and $h$ be integers with $c \ge 1$, $d \ge 2$ and 
$0 \le h \le h_{c,d,d} = \frac{(c+d-2)\cdots (c+1)c}{d!}$. 
Also, let $\Delta^{\rm{min}}$ be a $(d-1)$-dimensional Buchsbaum 
simplicial complex with $d$-linear resolution and 
$\dim_k H_{\frm}^{d-1}(k[\Delta^{\rm{min}}]) 
= \lfloor h_{c,d,d}\rfloor$. 
Then we can prove the following$:$

\par \vspace{1mm} \par \noindent 
{\bf Corollary \ref{MinCover}.} 
If such a complex $\Delta^{\rm{min}}$ exists, then  there exists a 
Cohen--Macaulay cover $\widetilde{\Delta}$ of $\Delta^{\rm{min}}$.  

\par \vspace{1mm} \par \noindent 
{\bf Theorem \ref{Int-d-linear}.}
Let $\Delta^{-} \subseteq \Delta \subseteq \Delta^{+}$ be simplicial 
complexes on $V=\{x_1,\ldots,x_n\}$. 
If both $k[\Delta^{-}]$ and $k[\Delta^{+}]$ are Buchsbaum 
Stanley--Reisner rings with $d$-linear  
resolutions, then so is $k[\Delta]$.  

\par \vspace{1mm} 
Thus we can reduce Conjecture \ref{Conj-1} to the existence 
of $\Delta^{\rm{min}}$. 
For $d=3$, we show that Hanano's example 
is turned out to be $\Delta^{\rm{min}}$ (Example \ref{Hanano-Ex}).  
Combining the above result with Corollary \ref{Conj-2linear}
and Example \ref{Cyclic-MM}, we have$:$

\par \vspace{1mm} \par \noindent 
{\bf Theorem.} 
Conjecture \ref{Conj-1} is true if either one of the following 
conditions are satisfied$:$ 
\par \vspace{1mm}
(1) $d \le 3$. \qquad (2) $q=2$. \qquad (3) $h_{c,d,q} =1$. 
 
\section{Preliminaries}

\par 
We first fix notation. 
Let $\N$ (resp. $\Z$) denote the set of nonnegative integers 
(resp. integers). 
Put $[n] = \{1,2,\ldots,n\}$ for any positive integer $n$.   
Let $\#(W)$ denote the cardinality of a set $W$. 

\par
We recall some notation on simplicial complexes and Stanley--Reisner 
rings according to \cite{St}. 
We refer the reader to e.g. \cite{BrHe}, 
\cite{Hi1}, \cite{Hoc} and \cite{St} for the detailed information about 
combinatorial and algebraic background. 

\par \vspace{1mm}
\begin{nota}[\textbf{Simplicial complex, face}]
A \textit{simplicial complex} $\Delta$ on the \textit{vertex set}
 $V = \{x_1,\ldots,x_n\}$ is 
a collection of subsets of $V$ such that 
\begin{enumerate}
\item[(i)]$\{x_i\} \in \Delta$ for every $1 \le i \le n$; 
\item[(ii)] $F \in \Delta$, 
$G \subseteq F \Longrightarrow G \in \Delta$. 
\end{enumerate}
Each element $F$ of $\Delta$ is called a \textit{face} of $\Delta$. 
A face $F$ is called an \textit{$i$-face} if $\#(F) = i+1$. 
We set $d = \max\{\#(F) \,|\, \sigma \in \Delta\}$ and 
define the \textit{dimension} of $\Delta$ to be $\dim \Delta =d-1$.
A maximal face is called a \textit{facet}.  
We say that $\Delta$ is \textit{pure} if every facet has 
the same cardinality.   
\par 
Given a subset $W$ of $V$, the \textit{restriction} of $\Delta$ is the 
subcomplex 
\[
 \Delta_{W} = \{F \in \Delta \,|\, F \subseteq W\}. 
\] 
On the other hand, if $F$ is a face of $\Delta$, then we define 
the subcomplex $\link_{\Delta}(F)$, 
which is called the \textit{link of $F$} as follows$:$
\[
 \link_{\Delta} (F) = \{G \in \Delta \,|\, 
 F \cap G = \emptyset,\; F \cup G \in \Delta \}. 
\]
In particular, $\link_{\Delta}(\emptyset) = \Delta$. 
\end{nota}

\par \vspace{1mm}
In the following, let $\Delta$ be a $(d-1)$-dimensional 
simplicial complex on $V = \{x_1,\ldots,x_n\}$. 
Let $S = k[x_1,\ldots,x_n]$ be a polynomial ring 
in $n$-variables over a field $k$. 
Here, we identify each $x_i \in V$ with the indeterminate $x_i$ of $S$. 
Also, let $\frm_S = (x_1,\ldots,x_n)S$ denote 
the homogeneous maximal ideal of $S$.

\par \vspace{1mm}
\begin{nota}[\textbf{$h$-vector, Reduced homology}]
Let $f_i = f_i(\Delta)$, $0 \le i \le d-1$, denote the numbers of $i$-faces 
in $\Delta$. 
We define $f_{-1} = 1$. 
We call $f(\Delta) = (f_0,f_1,\ldots, f_{d-1})$ 
the \textit{$f$-vector} of $\Delta$. 
Also, we define the \textit{$h$-vector} 
$h(\Delta) = (h_0,h_1,\ldots, h_d)$ of $\Delta$ by 
\begin{equation}
 \sum_{i=0}^d  \frac{f_{i-1}t^i}{(1-t)^{i}} = 
 \frac{h_0 + h_1t + \cdots + h_{d}t^d}{(1-t)^d}.    
\end{equation}  
\par 
Let $C_p(\Delta)$ be the $k$-vector space generated by 
all $p$-faces of $\Delta$ and we define the differential map 
$\partial_{p+1} \colon C_{p+1}(\Delta) \to C_p(\Delta)$ as follows$:$
\[
  \partial_{p+1}(\{x_{i_1},\ldots, x_{i_{p+1}}\}) = 
\sum_{j=1}^{p+1} (-1)^j \{x_{i_1},\ldots, \widehat{x_{i_j}}, 
\ldots, x_{i_{p+1}}\}. 
\]
Then $C_{\bullet}(\Delta)$ is a bounded complex and 
$\widetilde{H}_p(\Delta;k) = H_{p}(C_{\bullet}(\Delta))$ is called 
the \textit{$i$th reduced simplicial homology group} of $\Delta$ 
with values in $k$. 
Note that $\widetilde{H}_{-1}(\Delta;k) = 0$ 
if $\Delta \ne \{\emptyset\}$ and 
\[
 \widetilde{H}_{p}(\{\emptyset\};k) =\left\{
 \begin{array}{cl}
 0  &  (p \ge 0); \\
 k  & (p = -1). 
 \end{array}
 \right.
\]
Also, note that $\widetilde{H}_p(\Delta;k) = 0$ 
if $p \ge d = \dim \Delta +1$ or $p \le -2$. 
\end{nota}

\par \vspace{1mm}
\begin{nota}[\textbf{Stanley--Reisner ring}]
Let $I_{\Delta}$ be the ideal of $S$ which is generated by 
square-free monomials $x_{i_1}\cdots x_{i_r}$, 
$1 \le i_1 < \cdots < i_r \le n$, with 
$\{x_{i_1},\cdots, x_{i_r}\} \notin \Delta$. 
We say that the quotient algebra $k[\Delta]:=S/I_{\Delta}$ is the 
\textit{Stanley--Reisner ring} of $\Delta$ over $k$. 
We consider $k[\Delta]$ as a graded algebra 
$k[\Delta] = \oplus_{j \ge 0} k[\Delta]_{j}$ with the standard grading, 
i.e., each $\deg x_i =1$.  
\par 
Let $e(k[\Delta])$ denote the \textit{multiplicity} of $k[\Delta]$. 
It is well known that 
$e(k[\Delta]) = f_{d-1}(\Delta) = \sum_{i=0}^{d} h_i(\Delta)$. 
\end{nota}

\begin{thm}[\textbf{Hochster's formula on the local cohomology modules}] 
\label{HocLoc}
Let $k[\Delta]$ be a Stanley--Reisner ring and $\frm$ 
the unique homogeneous maximal ideal of $k[\Delta]$. 
Then the Hilbert series of 
the $i$th local cohomology module $H_{\frm}^i(k[\Delta])$  
is given by 
\begin{equation}
 F(H_{\frm}^i(k[\Delta]),t) 
 = \sum_{F \in \Delta} \dim_k \widetilde{H}_{i-\#(F)-1}(\link_{\Delta} F; k) 
 \left(\frac{t^{-1}}{1-t^{-1}}\right)^{\#(F)}, 
\end{equation}
where $F(M,t) = \sum_{n \in \bbZ} \dim_k M_n t^n$ for a graded $S$-module 
$M$ with $\dim_k M_n < \infty$ for all $n \in \bbZ$.  
\par
In particular, we have 
\begin{enumerate}
 \item  $[H_{\frm}^i(k[\Delta])]_j = 0$ for all $j \ge 1$. 
 \item  $[H_{\frm}^i(k[\Delta])]_0 \cong \widetilde{H}_{i-1}(\Delta;k)$.  
 \item For all $p \ge 1$, 
\[
 \dim_k [H_{\frm}^i(k[\Delta])]_{-p} 
= \sum_{F \in \Delta,1 \le \#(F) \le p}
 \dim_k  \widetilde{H}_{i-\#(F)-1}(\link_{\Delta} F; k). 
\]
\end{enumerate}
\end{thm}

\par \vspace{1mm}
Also, if we define the \textit{$a$-invariant} 
of a $d$-dimensional graded ring $A$ 
(see \cite{GoWa}) by 
\[
 a(A) = \max\{m \in \Z \,|\, [H_{\frm}^d(A)]_m \ne 0  \}, 
\]
then we always have $a(k[\Delta]) \le 0$. 

\par \vspace{1mm}
\begin{nota}[\textbf{Cohen--Macaulay, Buchsbaum complex}]
$\Delta$ is called \textit{Cohen--Macaulay} over $k$ if 
it satisfies one of the following equivalent conditions$:$ 
\begin{enumerate}
 \item $k[\Delta]$ is Cohen--Macaulay. 
 \item  $H_{\frm}^i(k[\Delta]) = 0$ for all $i < d$. 
 That is, $\depth\,k[\Delta] = \dim k[\Delta]$. 
  \item $\widetilde{H}_i(\link_{\Delta}(F);k)=0$ 
   for every $F \in \Delta$ and for 
   every $i < \dim \link_{\Delta}(F)$. 
\end{enumerate}

\par \vspace{1mm}
Also, $\Delta$ is called \textit{Buchsbaum} over $k$ if 
it satisfies one of the following equivalent conditions$:$ 
\begin{enumerate}
 \item $k[\Delta]$ is Buchsbaum. 
 That is, $l(k[\Delta]/J) - e(J)$ (this invariant is 
 written as $I(k[\Delta])$) is independent on 
 the choice of homogeneous parameter ideal $J$ of $k[\Delta]$, 
 where $e(J)$ denotes the multiplicity of $J$.   
 \item $k[\Delta]$ is (\textit{F.L.C.}), i.e., 
  $l(H_{\frm}^i(k[\Delta])) < \infty$ for all $i < d$. 
 \item $H_{\frm}^i(k[\Delta]) = [H_{\frm}^i(k[\Delta])]_0 
  = \widetilde{H}_{i-1}(\Delta;k)$ for all $i < d$. 
 \item 
  $\Delta$ is pure and $\widetilde{H}_i(\link_{\Delta}(F);k)=0$ 
  for every $F (\ne \emptyset) \in \Delta$ and for 
  every $i < \dim \link_{\Delta}(F)$.  
\end{enumerate} 
When this is the case, 
\[
 I(k[\Delta]) = \sum_{i=0}^{d-1} \bbinom{d-1}{i} \,l(H_{\frm}^i(k[\Delta]))
 = \sum_{i=0}^{d-1} \bbinom{d-1}{i} \,\dim_k \widetilde{H}_{i-1}(\Delta;k).
\]
\par 
Note that $k[\Delta]$ is Cohen--Macaulay 
if and only if it is Buchsbaum and 
$\widetilde{H}_i(\Delta;k) =0$ for all $i < d-1$.  
\end{nota}

\par \vspace{1mm}
\begin{nota}[\textbf{Regularity, Linear resolution}]
Let $A = k[A_1] =S/I$ be a homogeneous $k$-algebra. 
\par
A graded \textit{minimal free resolution} (abbr., MFR) 
of $A$ over $S$ is an exact sequence 
\[
 0 \to \bigoplus_{j \in \Z} S(-j)^{\beta_{p,j}(A)} 
\stackrel{\varphi_{p}}{\longrightarrow} \cdots 
\stackrel{\varphi_{2}}{\longrightarrow}  
\bigoplus_{j \in \Z} S(-j)^{\beta_{1,j}(A)}  
\stackrel{\varphi_{1}}{\longrightarrow}  
S \to A \to 0,  
\]
where $S(j)$, $j \in \Z$, denotes the graded 
module $S(j) = \oplus_{n \in \Z} S_{j+n}$ and \lq\lq minimal'' means 
that $\varphi_i \otimes_A A/\frm = 0$ for all $i$. 
We say that $\beta_i(A) = \sum_{j \in \Z} \beta_{i,j}(A)$ 
the $i$th \textit{Betti number} of $A$ over $S$. 

\par
The \textit{Castelnuovo--Mumford regularity} is defined by 
\[
 \reg A := \max\{j-i \,|\, \beta_{i,j}(A) \ne 0\}. 
\]
Also, the \textit{initial degree} of $A$ is defined by
\[
 \indeg A := \min\{j \,|\, I_j \ne 0\} 
= \min\{j \,|\, \beta_{1,j}(A) \ne 0 \}.
\]
\par 
Notice the following fact$:$
\begin{enumerate}
 \item $\beta_{i,j}(A) = \dim_k \Tor_S^i(A,S/\frm_S)_j$. 
 \item $p = n -\depth A$ by Auslander--Buchsbaum formula.  
 \item $\beta_{i,j}(A) = 0$ for all $j < i+\indeg A-1$. 
 \item $\reg A \ge \indeg A -1$. 
 \item $\reg A  = \inf\{r \in \Z\,|\, 
\mbox{$[H_{\frm}^i(A)]_j =0$ for all $j > r-i$} \}$; see \cite{EiGo}. 
\end{enumerate}

\par \vspace{1mm} 
$A$ has \textit{$q$-linear resolution} (abbr., $A$ is $q$-linear or 
$\Delta$ is $q$-linear over $k$)  
if $\reg A = \indeg A -1 = q-1$ (e.g., \cite{Oo}), that is, 
its graded minimal free resolution of $A$ is written as 
the following form$:$
\[
  0 \to S(-(q+p-1))^{\beta_{p}} \to S(-(q+p-2))^{\beta_{p-1}} \to \cdots 
  \to S(-q)^{\beta_{1}} \to S \to A \to 0.  
\]
\end{nota}

\par
In the case of Buchsbaum homogeneous $k$-algebras, the following criterion
for having a $q$-linear resolution is known. 
See also Theorems \ref{HibiCri} and \ref{AinLink}. 

\begin{thm}[{\rm cf. \cite[Corollary 1.5]{EiGo}}] \label{MinRed}
Let $A=k[A_1]$ be a homogeneous Buchsbaum $k$-algebra 
with $\indeg A \ge q$
and $\dim A =d$. Put $\frm = A_{+}$, the unique homogeneous 
maximal ideal of $A$. 
Suppose that $k$ is infinite. 
Then the following conditions are equivalent$:$ 
\begin{enumerate}
 \item $A$ has $q$-linear resolution. 
 \item There exists a homogeneous system of parameters $f_1,\ldots,f_d$ 
   such that $\frm^q = (f_1,\ldots,f_d)\frm^{q-1}$. 
 \item $[H_{\frm}^i(A)]_j = 0$ for all $i \ne d$, $j \ne q-1-i$  
  and $[H_{\frm}^d(A)]_j = 0$ for all $j > q-d-1$.   
\end{enumerate}
\end{thm}

\par \vspace{1mm}
The Betti numbers of Stanley--Reisner rings
can be calculated by the following formula. 

\begin{thm}[\textbf{Hochster's formula on the Betti numbers} 
{\rm \cite[Theorem 5.1]{Hoc}}] \label{HochsterBetti}
Let $k[\Delta]$ be the Stanley--Reisner ring of $\Delta$ 
over $k$. 
Then 
\begin{equation}
 \beta_{i,j}(k[\Delta]) = 
\sum_{\begin{subarray}{c} W \subseteq V \\
\#(W) = j\end{subarray}} 
\!\!\!\! \dim_k \widetilde{H}_{j-i-1}(\Delta_W;k).  
\end{equation}
\end{thm}

\par \vspace{1mm}
Also, we frequently use the following theorem; 
see \cite[Corollary 2.8]{HoMi}. 

\begin{thm}[\textbf{Hoa--Miyazaki theorem}] \label{HoaMiy}
Let $A=k[A_1]$ be a $d$-dimensional homogeneous Buchsbaum $k$-algebra. 
Then
\[
 \reg A \le a(A) + d+1. 
\] 
\end{thm}
\vspace{2mm}
\section{Buchsbaum Stanley--Reisner rings with $q$-linear resolution}

\par
In this section, let us gather several properties of 
Buchsbaum Stanley--Reisner rings having a $q$-linear resolution.  
\par
Throughout this section, let $\Delta$ be a simplicial complex 
on $V = \{x_1,\ldots,x_n\}$, and 
let $k[\Delta]$ denote the Stanley--Reisner ring of $\Delta$ over a field $k$, 
and $\frm$ its homogeneous maximal ideal. 
Also, let $c$, $d$, $q$ be given integers with 
$c \ge 1$ and $q,d \ge 2$. 
\par 
Let us start giving fundamental results on Buchsbaum Stanley--Reisner rings 
with linear resolution. 

\begin{prop} {\rm (cf. \cite{EiGo})} \label{Known} 
Suppose that $A = k[\Delta]$ is a $d$-dimensional 
Buchsbaum Stanley--Reisner ring 
with $q$-linear resolution, and put $\codim A=c$. 
Then  
\begin{enumerate}
 \item $q \le d+1$. 
 \item $H_{\frm}^{i}(A) = [H_{\frm}^{i}(A)]_0 
 \cong \widetilde{H}_{i-1}(\Delta;k) =0$ 
  for all $i \ne q-1,\,d$. 
 \item $\reg A = q-1$ and $\;\indeg A = q$. 
 \item $a(A) = q-d-2$ or $q-d-1$. 
\end{enumerate}
\par 
Also, suppose that  $q \le d$. Then 
\begin{enumerate}
 \item[(5)] $\widetilde{H}_{d-1}(\Delta;k) =0$. 
 \item[(6)] $\dim_k H_{\frm}^{q-1}(A) 
  = \dim_k \widetilde{H}_{q-2}(\Delta;k) 
  =\beta_{c+d-q+1}$. 
\end{enumerate}
\end{prop}

\begin{proof}[\quad Proof]
(3) is clear by definition. 
By Hochster's formula, we have $\reg A \le d$ in general. 
Thus we get (1). 
Also, (2) easily follows from Theorem \ref{MinRed}. 
Similarly we have $a(A) \le q-d-1$. 
On the other hand, $q-1 = \reg(A) \le a(A) + d+1$ by Hoa--Miyazaki Theorem.  
This implies that $a(A) \ge q-d-2$. 
Hence we get (4). 
\par 
Now, suppose that $q \le d$. 
Since $\reg A = d > q-1$, we have that $\widetilde{H}_{d-1}(\Delta;k) 
= [H_{\frm}^d(A)]_0 = 0$. 
Also, (6) follows from Theorem \ref{HochsterBetti}
if we put $j = c+d=n$ and $i=c+d-q+1$. 
\end{proof}

\par
In the following, we always assume that $q \le d$. 
In fact, we can completely characterize Stanley--Reisner rings 
with $(d+1)$-linear resolution as follows$:$

\begin{prop} \label{d+1-linear}
Let $A = k[\Delta]$ be a $d$-dimensional Stanley--Reisner ring on $V$.  
Then the following conditions are equivalent$:$ 
\begin{enumerate}
 \item $A$ has $(d+1)$-linear resolution. 
 \item $\indeg A = d+1$. 
 \item $I_{\Delta} = (x_{i_1}\cdots x_{i_{d+1}}\,|\, 
  1 \le i_1 < \cdots < i_d < i_{d+1} \le n)$. 
  That is, $\Delta$ is the $(d-1)$-skeleton of the  
  standard $n$-simplex $2^{V}$.
\end{enumerate}
When this is the case, $A$ is Cohen--Macaulay. 
\end{prop}

\begin{proof}[\quad Proof]
$(1) \Longrightarrow (2)$ is trivial. 
Conversely, suppose (2). 
Then $d \ge \reg A \ge \indeg A -1 = d$. 
Thus $A$ has $(d+1)$-linear resolution. 
\par 
$(3) \Longrightarrow (2)$ is trivial. 
Conversely, suppose that $\indeg A = d+1$.
Then $I_{\Delta}$ does not contain any square-free 
monomial $M$ with $\deg M \le d$. 
Now suppose that some square-free monomial $x_{i_1}\cdots x_{i_{d+1}}$ 
of degree $(d+1)$ is not contained in $I_{\Delta}$.  
Then $\{x_{i_1},\ldots,x_{i_{d+1}}\} \in \Delta$. 
This contradicts the assumption that $\dim \Delta = d-1$. 
Hence $\Delta$ is a $(d-1)$-skeleton of $2^{V}$. 
Since $2^{V}$ is Cohen--Macaulay, so is $\Delta$; 
see e.g., \cite[Ex 5.1.23]{BrHe}. 
\end{proof}

\par
The following theorem gives a criterion for having a
$q$-linear resolution. 
In Section 4, we will improve this theorem. 

\begin{thm}[\textbf{Hibi's criterion} {\rm \cite[Theorem (1.6)]{Hi2}}] 
\label{HibiCri}
Let $A = k[\Delta]$ be a $d$-dimensional Buchsbaum Stanley--Reisner ring. 
Put $\indeg A = q \le d$. 
Then $A$ has $q$-linear resolution if and only if the following 
conditions are satisfied$:$
\begin{enumerate}
 \item[(i)] $\widetilde{H}_{i}(\Delta;k) = 0$ for all $i \ne q-2$. 
 \item[(ii)] $a(k[\link_{\Delta}(\{x_i\})]) \le q-d$ for all $i=1,\ldots,n$. 
\end{enumerate}
\end{thm}

\par \vspace{1mm}
Now suppose that $A = k[\Delta]$ is a $d$-dimensional Buchsbaum 
Stanley--Reisner ring with $q$-linear resolution. 
Then $h := \dim_k \widetilde{H}_{q-2}(\Delta;k) 
= \dim_k H^{q-1}_{\frm}(A)$ is an important 
invariant of $\Delta$.   
From now on, we focus this invariant. 
\par 
The following proposition may be known, but we give a proof 
for the readers' convenience. 

\begin{prop} \label{LinMul}
Suppose that $A = k[\Delta]$ is a $d$-dimensional 
Buchsbaum Stanley--Reisner ring with 
$q$-linear resolution, and put $\codim A=c$. 
If we put $h := \dim_k H^{q-1}_{\frm}(A)$, 
then the multiplicity $e(A)$ and the $I$-invariant $I(A)$ are given by 
\begin{equation}
 e(A)= \bbinom{c+q-1}{q-1} - h \bbinom{d-1}{q-1}  \quad \text{and}
 \quad I(A) = h \bbinom{d-1}{q-1}. 
\end{equation}
\end{prop}

\par \vspace{1mm}
To prove the above proposition, we first show the following lemma. 
See also \cite{EiGo}. 

\begin{lemma} \label{MulIneq}  
Let $A = k[A_1]$ be a $d$-dimensional homogeneous Buchsbaum $k$-algebra with 
$\codim A = c$ and $\indeg A \ge q$. 
Then 
\begin{enumerate}
 \item $e(A) \ge \sbinom{c+q-1}{q-1} - I(A)$. 
 \item Equality holds in $(1)$ if and only if 
 $A$ has $q$-linear resolution. 
\end{enumerate}
\end{lemma}

\begin{proof}[\quad Proof]
We may assume that $k$ is infinite. 
\par
(1) Let $J$ be a homogeneous minimal reduction of $\frm$, that is, 
$J$ is a homogeneous parameter ideal of $A$ and $\frm^{r+1} = J\frm^r$ holds 
for some integer $r \ge 0$. 
Then $e(A) = e(J) = l_A(A/J) - I(A)$ since $A$ is Buchsbaum. 
Also, since 
$B:=A/J$ is a homogeneous Artinian $k$-algebra with $\dim_k B_1 = c$ and 
$\indeg B \ge q$, we have 
\[
 e(A) = l_A(A/J) - I(A) \ge l_A(A/J+\frm^q) - I(A) 
 = \sbinom{c+q-1}{q-1} - I(A). 
\]
\par 
(2) Equality holds in $(1)$ if and only if $\frm^q \subseteq J$.  
Since $A$ is homogeneous, this yields that 
$\frm^q = J \cap \frm^q = J\frm^{q-1}$. 
Thus the statement follows from  Theorem \ref{MinRed}. 
\end{proof}

\begin{proof}[\quad Proof of Proposition $\ref{LinMul}$]
Since $A = k[\Delta]$ has $q$-linear resolution, 
$H_{\frm}^i(A) = 0$ for all $i \ne q-1,d$ and thus 
$I(A) = h \bbinom{d-1}{q-1}$. 
On the other hand, by Lemma \ref{MulIneq}, we have  
\[
  e(A) = \bbinom{c+q-1}{q-1} -I(A) 
  =  \bbinom{c+q-1}{q-1} - h \bbinom{d-1}{q-1}, 
\] 
as required. 
\end{proof}

\par
The following theorem plays an important role in this article.  

\begin{thm} \label{Key}  
Suppose that $A = k[\Delta]$ is a $d$-dimensional 
Buchsbaum Stanley--Reisner ring with 
$q$-linear resolution $(q \le d)$, and put $c=\codim A$ and 
$h = \dim_k H^{q-1}_{\frm}(A)$. 
Then the $h$-vector 
$(h_0,h_1,\ldots,h_{q-1},h_q,h_{q+1},\ldots,h_d)$ of $\Delta$ is 
\begin{equation}
  \left(1,\,c,\,\cdots, \sbinom{c+q-2}{q-1},\,- \sbinom{d}{q}h,
  \,\sbinom{d}{q+1}h,\cdots,(-1)^{d-q+1} \sbinom{d}{d}h\right). 
\end{equation}
That is, 
\[
 h_p = \bbinom{c+p-1}{p} \;\;(1 \le p \le q-1), \quad 
 h_p=  (-1)^{p-q+1}\bbinom{d}{p} h \;\;(q \le p \le d).
\]
\par 
Also, $h$ satisfies the following inequalities$:$
\begin{equation}
 0 \qle h \qle
 \frac{(c+q-2)\cdots (c+1)c}{d(d-1)\cdots (d-q+2)}=:h_{c,d,q}.
\end{equation}
\end{thm}

\begin{proof}[\quad Proof]
Since $\indeg A = q$,  one has 
\[
  h_p = \bbinom{c+p-1}{p} \quad \text{for all $p = 0,1,\ldots,q-1$}.  
\]
\par 
On the other hand, by the similar argument as in the proof of 
\cite[Theorem 2.1]{Te1}, we have 
\begin{align*}
 \dim_{k} [H_{\frm}^d(A)]_{-1} 
 & =  d \cdot h_d + h_{d-1} \\ 
 \dim_{k} [H_{\frm}^d(A)]_{-2} 
 & =  \bbinom{d+1}{2} h_d + d \cdot h_{d-1} + h_{d-2} \\ 
 \dots &  \dots \\
 \dim_k [H_{\frm}^d(A)]_{-p} 
 & =  \bbinom{d+p-1}{p}h_d + \bbinom{d+p-2}{p-1}h_{d-1} + \cdots 
 + d \cdot h_{d-p+1} + h_{d-p} \\ 
 \dots &  \dots \\
 \dim_k [H_{\frm}^d(A)]_{q-d} 
 & =  \bbinom{2d-q-1}{d-q}h_d + \bbinom{2d-q-2}{d-q-1}h_{d-1} 
 + \cdots + d \cdot h_{q+1} + h_{q}. 
\end{align*}
By Proposition \ref{Known}, we have 
\[
 (-1)^{d-1} h_d  = \widetilde{\chi}(\Delta) = \sum_{i=-1}^{d-1} (-1)^i 
 \dim_k \widetilde{H}_i(\Delta;k) = (-1)^q h
\]
and $\dim_k [H_{\frm}^d(A)]_{j} =0$ for all $j=-1,-2,\ldots,q-d$. 
Solving the above equations, one can easily 
obtain that $h_p =(-1)^{p-q+1} \bbinom{d}{p} h$ 
for all $p=q,\ldots,d-1,d$.  
Then 
\[
 \bbinom{2d-q}{d-q+1}h_d + \cdots + d \cdot h_{q} 
 + h_{q-1} 
 = \dim_k [H_{\frm}^d(A)]_{q-d-1} \ge 0
\]
implies that 
\[ 
  \bbinom{c+q-2}{q-1} - \bbinom{d}{q-1}h \ge 0.  
\]
Namely, $h \le h_{c,d,q}$, as required. 
\end{proof}

\begin{remark} \label{h-vec}
Proposition $\ref{LinMul}$ also follows from Theorem $\ref{Key}$. 
\end{remark}

\begin{remark}
Under the same notation as in Theorem $\ref{Key}$, 
let $\beta_i$ be the $i$th Betti number of $k[\Delta]$ over $S$. 
Then the following identity holds$:$ 
\[
 (1-t)^c \left\{\sum_{p=0}^{q-1} \bbinom{c+p-1}{p} t^p 
+ \sum_{p=q}^d (-1)^{p-q+1} \bbinom{d}{p}h t^p \right\}
= 1 - \sum_{i=1}^{n-q+1} \beta_i t^{q+i-1}. 
\]
\end{remark}

\par \vspace{1mm}
Based on Theorem \ref{Key}, we pose the following conjecture. 

\begin{conj} \label{Conj-1}
Let $d$, $c$, $q$, $h$ be integers with $c \ge 1$, $h \ge 0$, 
and $2 \le q \le d$. 
Then the following conditions  are equivalent$:$
\begin{enumerate}
\item There exists a Buchsbaum Stanley--Reisner ring $A=k[\Delta]$ 
with $q$-linear resolution
  such that  $\dim A =d$, $\codim A = c$ and 
  $\dim H_{\frm}^{q-1}(A)=h$.
\item The following inequality holds$:$ 
\[
  0 \le h \le 
 h_{c,d,q} = \frac{(c+q-2)\cdots (c+1)c}{d(d-1)\cdots (d-q+2)}.
\]
\end{enumerate}
\end{conj}

\begin{remark}[{\rm Stanley}] \label{CM-conj}
There is a $(d-1)$-dimensional Cohen--Macaulay complex $\Delta$ with 
h-vector 
\[
 (h_0,\ldots,h_{q-1}) = \left(1,c,\sbinom{c+1}{2}, \cdots, 
\sbinom{h+q-2}{q-1}\right), 
\]
which has $q$-linear resolution. See \cite[Theorem 5.1.15]{BrHe}. 
\par
In particular, Conjecture $\ref{Conj-1}$ is true in the case of $h=0$. 
\end{remark}

\par
Now let $A=k[A_1]$ be a $d$-dimensional Buchsbaum 
homogeneous $k$-algebra.  
Let $e(A)$ (resp. $\emb(A) = \dim_k A_1$) denote the multiplicity 
(resp. the embedding dimension) of $A$. 
Then $\emb(A) \le e(A) + \dim A + I(A)-1$ holds in general, 
and $A$ is said to \textit{have maximal embedding dimension} 
if equality holds. 
Also, $A$ has maximal embedding dimension 
if and only if it has $2$-linear resolution or is isomorphic to 
a polynomial ring; see e.g. \cite{Go1}. 

\par
Fr\"oberg (\cite{Fr1,Fr2}) has determined the structure 
of Buchsbaum  simplicial complexes with  
$2$-linear resolution. 

\begin{prop}[{\rm Fr\"oberg \cite{Fr2}}] \label{2-linearBbm}
Let $\Delta$ be a $(d-1)$-dimensional simplicial complex which is 
not a $(d-1)$-simplex.   
Then the following conditions are equivalent$:$ 
\begin{enumerate}
 \item $k[\Delta]$ is Buchsbaum with $2$-linear resolution. 
 \item $\Delta$ is a finite disjoint union of $(d-1)$-dimensional 
  simplicial complexes $\Delta_i$ such that 
 $k[\Delta_i]$ is Cohen--Macaulay of maximal embedding dimension.   
\end{enumerate}
\end{prop}

\par
Using the above proposition, we can show that the above conjecture 
is true for $q=2$. 

\begin{cor} \label{Conj-2linear}
Let $d$, $c$, $h$ be integers with $d \ge 2$, 
$c \ge 1$, and $h \ge 0$.
Then the following conditions  are equivalent$:$
\begin{enumerate}
\item There exists a Buchsbaum Stanley--Reisner ring $A=k[\Delta]$ 
  with $2$-linear resolution
  such that  $\dim A =d$, $\codim A = c$ and 
  $\dim H_{\frm}^{1}(A)=h$.
\item The following inequality holds$:$ 
\[
  0 \le h \le h_{c,d,2} = \frac{c}{~d~}.
\]
\end{enumerate}
\end{cor}

\begin{proof}[\quad Proof]
Note that there are many examples of Cohen--Macaulay 
Stanley--Reisner ring $A$ over any field $k$ 
of maximal embedding dimension such 
that $\codim A = c$ and $\dim A =d$. 
Indeed, $k[X_0,\ldots,X_c,Y_1,\ldots,Y_{d-1}]/(X_iX_j\,|\,0 \le i < j \le c)$
gives one of such examples. 
\par
To see $(2)\Longrightarrow (1)$, suppose that $c \ge dh$. 
We may assume that $h > 0$. 
Take any simplicial complex $\Delta_0$ for which 
$k[\Delta_0]$ is Cohen--Macaulay of maximal embedding dimension 
and $\codim k[\Delta_0] =c-dh$.
Let $\Delta$ be a disjoint union of $\Delta_0$ and $(d-1)$-simplexes 
$\Delta_1,\ldots,\Delta_h$. 
Then $A=k[\Delta]$ is Buchsbaum with $2$-linear resolution 
by Proposition \ref{2-linearBbm}. 
Also, we have $\codim A = (c-dh) + d \cdot h = c$,  
$\dim A = d$ and $\dim_k H_{\frm}^1(A) = h$. 
\end{proof}

\par 
On the other hand, in the case of $q \ge 3$, 
it seems to be difficult to construct examples 
of complexes having a $q$-linear resolution with given parameters.
But, in Section 5, we will give an affirmative answer in the 
case of $q=d=3$ using the notion 
of ``\textit{Cohen-Macaulay cover}''.

\vspace{2mm}
\section{Alexander Duality of Buchsbaum complex with linear resolution}

In this section, we characterize 
the Alexander dual of Buchsbaum simplicial complexes
with linear resolution. 
As an application, we give some examples of 
Buchsbaum complexes with linear resolution
using cyclic polytopes. 
We first recall some basic results on Alexander duality of 
simplicial complexes. 

\begin{nota}[\textbf{Alexander duality}] \label{Alex}
Let $\Delta$ be a $(d-1)$-dimensional simplicial complex on 
$V = \{x_1,\ldots,x_n\}$. 
The \textit{Alexander dual} of $\Delta$ is defined by 
\[
 \Delta^{*} := \{F \subseteq V \,|\, V \setminus F \notin \Delta \}. 
\]
Suppose that $c := n-d \ge 2$ and $\indeg A  = q \ge 2$. 
Then $\Delta^{*}$ is a simplicial complex on $V$. 
Also, the following statements hold$:$ 
see \cite{EaRe, Te2} for details. 
\begin{enumerate}
 \item (Alexander Duality) 
  $\widetilde{H}_{i-2}(\Delta^{*};k) \cong 
  \widetilde{H}^{n-i-1}(\Delta;k)$ for all $i$. 
 \item $(\Delta^{*})^{*} = \Delta$. 
 \item $\dim k[\Delta^{*}] + \indeg k[\Delta]=c+d =n$. 
 In particular, 
 \[
  \dim k[\Delta^{*}] = c+d-q,\quad \indeg k[\Delta^{*}] = c, 
  \quad \text{and}\quad \codim k[\Delta^{*}] = q. 
 \]
 \item The Betti numbers of $k[\Delta^{*}]$ are given by the formula
\[
 \beta_{i,j}(k[\Delta^{*}]) 
 = \sum_{\begin{subarray}{c} F \in \Delta \\ \#(F) = c+d-j \end{subarray}} 
 \dim_k \widetilde{H}_{i-2}(\link_{\Delta}F; k).
\]
 \item $k[\Delta]$ has linear resolution if and only if  
$k[\Delta^{*}]$ is Cohen--Macaulay.  
In fact, we have 
\[
  \reg k[\Delta] - \indeg k[\Delta]+1 
= \dim k[\Delta^{*}] - \depth k[\Delta^{*}]. 
\]
\end{enumerate}
\end{nota}

\par \vspace{1mm}
Described as above, if $\Delta$ is Cohen--Macaulay with 
linear resolution, then so is $\Delta^{*}$. 
Thus it is a natural to ask 
\par
\lq\lq What can we say about the Alexander dual 
of a \textit{Buchsbaum} simplicial complex
with linear resolution ?'' 
\par
An answer to this question is that 
``The Alexander dual of such a complex has  
\textit{almost linear resolution} (see \cite{EiGo}) 
with suitable conditions''. 
To be precise, we have$:$

\begin{thm} \label{LinAlex}
Let $c,d,q$ be integers with $c \ge 2$, $2 \le q \le d$. 
Let $A = k[\Delta]$ be a $d$-dimensional 
Stanley--Reisner ring with $\codim A = c$ and $\indeg A =q$,
and let $\Delta^{*}$ denote the Alexander dual of $\Delta$. 
Put $A^{*}:= k[\Delta^{*}]$. 
Then the following conditions are equivalent$:$ 
\begin{enumerate}
\item $A$ is Buchsbaum with $q$-linear resolution. 
\item $A^{*}$ is Cohen--Macaulay with 
 almost $c$-linear resolution and   
 the graded minimal free resolution of $A^{*}$ over $S$ can be written 
 as follows$:$
\[
 0 \to F_q \to F_{q-1} 
 = S(-(c+q-2))^{\beta_{q-1}^{*}} \to \cdots \to F_1 = S(-c)^{\beta_{1}^{*}} 
 \to S \to A^{*} \to 0,   
\]
where $F_{q} = S(-(c+d))^{\beta^{*}} \oplus S(-(c+q-1))^{\beta^{*'}}$. 
\end{enumerate}
When this is the case, $\beta^{*} = \dim_k H_{\frm}^{q-1}(A)$ and 
$\beta^{*'} = \dim_k H_{\frm}^d(A)_{q-d-1}$.  
\end{thm}

\begin{proof}[\quad Proof]
We may assume that $A$ has $q$-linear resolution and $A^{*}$ is 
Cohen--Macaulay by \ref{Alex}. 
\par
$(1)\Longrightarrow (2):$
To see (2), we must show that 
$\beta_{ij}(k[\Delta^{*}]) = 0$
for all pairs $(i,j)$ with $0 \le i \le q$ and $j > c+i-1$ 
except $(i,j) =(q,c+d)$. 

\par 
Let $F$ be a face of $\Delta$ with $\#(F) = c+d-j$. 
First suppose that $j=c+d$. 
Then $F=\emptyset$. 
If $i \le q-1$, then $\widetilde{H}_{i-2}(\Delta;k)=0$ by 
Proposition \ref{Known}.  
Next suppose that $c+i-1 < j < c+d$. 
Then $F \ne \emptyset$. 
Since $i-2 < j-c-1=d-\#(F)-1 = \dim \link_{\Delta}F$, 
the Buchsbaumness of $\Delta$ implies that
$\widetilde{H}_{i-2}(\link_{\Delta} F ;k)= 0$. 
Thus we get the required vanishing. 

\par 
$(2)\Longrightarrow (1):$
First, note that $\Delta$ is pure.  
Indeed, $I_{\Delta^{*}}$ is minimally generated by the elements 
$x_{j_1}\cdots x_{j_c}$ for which $V \setminus \{x_{j_1},\ldots,x_{j_c}\}$ 
is a maximal face of $\Delta$. 
\par
To see the Buchsbaumness of $k[\Delta]$, let 
$F$ be any non-empty face of $\Delta$ and let $i$ be 
an integer with $i < d-\#(F)+1$. 
Put $j = c+d-\#(F)$. 
Then one can easily see that $\beta_{ij}(k[\Delta^{*}]) =0$ 
by the assumption (2). 
This implies that $\widetilde{H}_{i-2}(\link_{\Delta}(F);k) =0$. 
Hence $k[\Delta]$ is Buchsbaum, as required.  
\par 
Putting $i=q$, $j=n(=c+d)$, we have 
\[
 \beta^{*} = \beta_{q,n}(k[\Delta^{*}]) 
= \dim_k \widetilde{H}_{q-2}(\Delta;k) 
= \dim_k H_{\frm}^{q-1}(A). 
\]
On the other hand, since $a(k[\Delta]) \le q-d-1$, 
we have $\widetilde{H}_{d-\#(F)-2}(\link_{\Delta}F; k)=0$
for every face $F$ with $1 \le \#(F) \le d-q$.    
Thus by Hochster's formula we have 
\begin{eqnarray*}
 \dim_k \left[H_{\frm}^d(k[\Delta])\right]_{q-d-1} 
\! & = & \!\!\!\!\!
\sum_{F \in \Delta,\#(F) = d-q+1} \!\!\!\!\!\!
\dim_k \widetilde{H}_{q-2}(\link_{\Delta}F; k) \\
&= & \beta_{q,c+q-1}(k[\Delta^{*}]) = \beta^{*'}. 
\end{eqnarray*}
\end{proof}

\begin{cor} \label{Alex-Gor}
Let $c$, $q$ and $d$ be integers with $c \ge 2$ and $2 \le q \le d$. 
Put $n = c+d$.  
Let $\Gamma$ be a simplicial complex on $V$. 
Suppose that $k[\Gamma]$ is a $(n-q)$-dimensional Gorenstein ring 
with almost $c$-linear resolution and $a(k[\Gamma]) =0$.  
Let $\Delta = \Gamma^{*}$ be the Alexander dual of $\Gamma$.  
Then $k[\Delta]$ is a $d$-dimensional Buchsbaum Stanley--Reisner ring 
with $q$-linear resolution. 
\end{cor}

\begin{proof}[\quad Proof]
Since $k[\Gamma]$ is Gorenstein and thus is Cohen--Macaulay, 
the length of the graded MFR of $k[\Gamma]$ is equal to 
$n-\dim k[\Gamma] = q$ by Auslander--Buchsbaum formula. 
Also,  since the graded MFR of $k[\Gamma]$ is almost $c$-linear, it 
can be written as follows$:$
\[
 0 \to F_q \to F_{q-1} 
 = S(-(c+q-2))^{\beta_{q-1}^{*}} \to \cdots \to F_1 = S(-c)^{\beta_{1}^{*}} 
 \to S \to k[\Gamma] \to 0,   
\]
where $S = k[X_1,\ldots,X_n]$ and $F_q = S(-\epsilon)$. 
Then $\epsilon - n = a(k[\Gamma]) = 0$; hence $\epsilon = n = c+d$. 
Hence $k[\Delta]$ has $q$-linear resolution by Theorem \ref{LinAlex}. 
\end{proof}

\par \vspace{1mm}
Let $n$, $f$ be integers with $n \ge f+1$. 
Consider the algebraic curve $M \subseteq \bbR^f$, defined 
by parametrically by $x(t) = (t,t^2,\ldots,t^{f})$, $t \in \bbR$. 
Let $C(n,f)$ be the convex hull of any distinct 
$n$-points over $M \subseteq \bbR^f$. 
Then $C(n,f)$ becomes a simplicial $f$-polytope. 
It is called  a {\it cyclic polytope} with $n$ vertices; 
\par 
It is well-known 
that any simplicial $f$-polytope $P$ with $n$ vertices 
satisfies $0 \le h_i(P) \le \bbinom{n-f+i-1}{i} = h_i(C(n,f))$; 
see e.g., {\rm \cite[Section 5.2]{BrHe}} for details.

\begin{exam}[\textbf{The Alexander dual of a cyclic polytope}] 
\label{CyclicPoly}
Let $q$, $d$ be integers with $2 \le q \le d$. 
Put $n = 2d-q+2$ and $f = 2(d-q+1)$. 
Also, let $\Gamma = \Gamma_{n,f}$ be the boundary complex of 
a cyclic polytope $C(n,f)$, 
and let $\Delta = \Gamma^{*}$ be the Alexander dual of $\Gamma$. 
Then $k[\Delta]$ is a $d$-dimensional Buchsbaum Stanley--Reisner ring 
with $q$-linear resolution with $h=1$. 
\end{exam}

\begin{proof}[\quad Proof]
Since $\Gamma$ is a boundary complex of a simplicial $f$-polytope $C(n,f)$, 
$k[\Gamma]$ is a $f$-dimensional Gorenstein Stanley--Reisner 
ring with $a(k[\Gamma]) = 0$. 
Indeed, $[H_{\frm}^d(k[\Gamma])]_0 
\cong \widetilde{H}_{d-1}(\Gamma;k) \cong k \ne 0$. 
Also, we have $\indeg k[\Gamma] = f/2+1 = d-q+2=n-d$ 
since $h_i(\Gamma) =\bbinom{n-f+i-1}{i}$ for all $i$.  
This implies that $k[\Gamma]$ has almost $(n-d)$-linear 
resolution (see \cite{Sch, TeHi}).  
On the other hand, we have 
\[
 \dim k[\Delta] = n - \indeg k[\Gamma] = d,\qquad 
 \indeg k[\Delta] = n - \dim k[\Gamma] = n-f = q. 
\]
Thus the assertion easily follows from the above corollary and 
Theorem \ref{LinAlex}. 
\end{proof}

\vspace{2mm}
\section{Buchsbaum Stanley--Reisner rings with minimal multiplicity}

\par
Let $A=k[A_1]$ be a homogeneous $k$-algebra of dimension $d$
with the unique homogeneous maximal ideal $\frm=A_{+}$. 
In \cite{Go2}, Goto proved an inequality 
\[
 e(A) \ge 1 + \sum_{i=1}^{d-1} \bbinom{d-1}{i-1} 
 l_A(H_{\frm}^i(A))
\]
and called the ring $A$ 
\textit{a Buchsbaum ring with minimal multiplicity} 
if equality  holds.  
Also, he proved that a Buchsbaum homogeneous $k$-algebra 
with minimal multiplicity has maximal embedding dimension, 
and hence has $2$-linear resolution. 

\par 
In this section, in the class of Stanley--Reisner rings, 
we introduce the notion of \textit{Buchsbaum ring with minimal 
multiplicity of type $q$} and prove that such a ring has 
a $q$-linear resolution; see Theorem \ref{Main}.  
Furthermore, we give several characterizations of this notion. 
In the following, 
let $c$, $d$, $q$ be integers with $c \ge 2$, $2 \le q \le d$. 

\begin{defn}[\textbf{Minimal multiplicity of type $q$}]  \label{MinMul}
Let $A:=k[\Delta]$ be a $d$-dimensional Buchsbaum Stanley--Reisner ring 
with $\codim A=c$ and $\indeg A = q$. 
Then we say that $A$ has {\it minimal multiplicity of type $q$} if 
\[
  e(A) = \frac{c+d}{d} \bbinom{c+q-2}{q-2}. 
\]
\end{defn}

\begin{prop} \label{GotoMinimal}
Let $A=k[\Delta]$ be a $d$-dimensional 
Buchsbaum Stanley--Reisner ring with $\indeg A \ge 2$. 
Then the following conditions are equivalent$:$
\begin{enumerate}
 \item $A$ has minimal multiplicity in the sense of Goto \cite{Go2}.
 \item $A$ has minimal multiplicity of type $2$.  
 \item $\Delta$ is a finite disjoint union of $(d-1)$-simplexes. 
\end{enumerate}
When this is the case, 
the normalization $B$ of $A$ is isomorphic to 
a finite product of polynomial rings with $d$-variables.   
In particular, the number of connected components of $\Delta$ is 
equal to the multiplicity of $A$.  
\end{prop}

\begin{proof}[\quad Proof]
$(1) \Longrightarrow (2):$ 
Suppose that $A$ has minimal multiplicity. 
Then since $A$ has $2$-linear resolution, 
$H_{\frm}^i(A) = 0$ for all $i \ne 1,d$ by Proposition \ref{Known}. 
Thus $e:=e(A) = 1+h$ and $I(A) = (d-1)h$, 
where $h = \dim_k H_{\frm}^1(A)$. 
Also, since $A$ has maximal embedding dimension, we have 
\[
 n = \emb(A) = e+d-1+I(A) = e+d-1+(d-1)(e-1) = de. 
\] 
Hence $e = \frac{n}{d} = \frac{c+d}{d}$, as required. 
\par \vspace{1mm}
$(2) \Longleftrightarrow (3):$ 
Note that $e$ is equal to the number of facets of $\Delta$. 
Since each facet of $\Delta$ is a $(d-1)$-simplex, we have 
$n \le de$ by counting the vertices of $\Delta$. 
Furthermore, equality holds if and only if 
$\Delta$ is a disjoint union of all facets. 
\par \vspace{1mm}
$(3) \Longrightarrow (1):$ 
Let $\Delta = \Delta_0 \cup \cdots \cup \Delta_h$ be a simplicial complex 
such that each $\Delta_i$ is a $(d-1)$-simplex. 
By Proposition \ref{2-linearBbm}, $A = k[\Delta]$ is a Buchsbaum 
ring of maximal embedding dimension. 
In particular, $H_{\frm}^i(A)=0$ for all $i\ne 1,d$,  
and $h = \dim_k \widetilde{H}_{0}(\Delta;k) = \dim_k H_{\frm}^1(A)$. 
Also, since $\Delta$ has $(1+h)$ facets, we get 
\[
 e(A) = 1+h = 1 + \sum_{i=1}^{d-1} \genfrac{(}{)}{0pt}{}{d-1}{i-1} l_A(H_{\frm}^i(A)).
\]
Hence $A$ has minimal multiplicity. 
\end{proof}

\par
The following theorem, which is a main theorem in this article, 
will justify our definition of minimal multiplicity of type $q$.  

\begin{thm} \label{Main}
Let $A := k[\Delta]$ be a Buchsbaum Stanley--Reisner ring 
with $\codim A = c$ and $\indeg A = q$. 
Then 
\begin{enumerate}
 \item The following inequality holds$:$
\begin{equation}
   e(A) \ge \frac{c+d}{d} \bbinom{c+q-2}{q-2}. 
\end{equation}
 \item If $A$ has minimal multiplicity of type $q$, 
 then it has $q$-linear resolution. 
\end{enumerate}
\end{thm}

\begin{proof}[\quad Proof]
(1) Let $V = \{x_1,\ldots,x_n\}$ be the vertex set of $\Delta$
where $n=c+d$. 
Put $\Gamma_i = :\link_{\Delta}(\{x_i\})$, and let 
$\frm_i$ be the homogeneous maximal ideal of $k[\Gamma_i]$ for all $i$. 
Then $k[\Gamma_i]$ is a $(d-1)$-dimensional Cohen--Macaulay 
ring since $A$ is Buchsbaum. 
Also, we have that  
$\codim k[\Gamma_i] = c$ and $\indeg k[\Gamma_i] \ge q-1$ by the assumption. 

By Lemma \ref{MulIneq} (or see \cite[Corollary 1.11]{EiGo}), we get 
\[
  e(k[\Gamma_i]) \ge \bbinom{c+(q-1)-1}{(q-1)-1} = 
  \bbinom{c+q-2}{q-2}. 
\]
On the other hand, counting the number of facets of $\Delta$, 
we have 
\[
 d \cdot e(A) = \sum_{i=1}^n e(k[\Gamma_i]) \ge (c+d) \bbinom{c+q-2}{q-2}, 
\]
as required. 
\par\vspace{1mm}
(2) Suppose that the equality holds. 
Then $e(k[\Gamma_i]) = \bbinom{c+q-2}{q-2}$ for all $i$.
It follows from Proposition \ref{LinMul} that 
$k[\Gamma_i]$ has $(q-1)$-linear resolution. 
This implies that $a(k[\Gamma_i]) = q-1-(d-1)-1 = q-d-1$. 
Also, we have that $\indeg A =q$ since 
$\indeg A \ge q$ and $\indeg k[\Gamma_i]=q-2$. 
Thus the assertion follows from Theorem \ref{AinLink} below. 
\end{proof}

\begin{thm} \label{AinLink}
Let $A = k[\Delta]$ be a $d$-dimensional 
Buchsbaum Stanley--Reisner ring of 
$\Delta$ on $V = \{x_1,\ldots,x_n\}$. 
Put $\indeg A = q$. 
\begin{enumerate}
 \item {\bf (Improved version of Hibi's criterion)}
 $A$ has $q$-linear resolution if and only if the following 
 conditions are satisfied$:$
 \begin{enumerate}
  \item $\widetilde{H}_{q-1}(\Delta;k) =0$. 
  \item $a(k[\link_{\Delta}\{x_i\}]) \le q-d$ for all $i = 1,\ldots,n$. 
 \end{enumerate} 
 \item If $a(k[\link_{\Delta}\{x_i\}]) = q-d-1$ for all $i$, 
 then $A$ has $q$-linear resolution and $a(A) = q-d-2$.
\end{enumerate} 
\end{thm}

\begin{proof}[\quad Proof]
Put $\Gamma_i = \link_{\Delta}(\{x_i\})$ for all $i = 1,\ldots,n$. 
Since $A = k[\Delta]$ is Buchsbaum, $k[\Gamma_i]$ is 
Cohen--Macaulay for all $i$ and 
$H_{\frm}^p(A) =[H_{\frm}^p(A)]_0 = \widetilde{H}_{p-1}(\Delta;k)$ for all 
$p \le d-1$. 
Also, by Hochster's formula, we have 
\begin{eqnarray*}
  F(H^d_{\frm}(A),t) 
  & = &  \sum_{F \in \Delta} 
        \dim_k \widetilde{H}_{d-\#(F)-1}(\link_{\Delta}F; k) 
  \left(\frac{t^{-1}}{1-t^{-1}}\right)^{\#(F)}; \\
  F(H^{d-1}_{\frm_i}(k[\Gamma_i]),t) 
  & = &  \sum_{G \in \Gamma_i} 
        \dim_k \widetilde{H}_{d-\#(G)-2}(\link_{\Gamma_i}G; k) 
  \left(\frac{t^{-1}}{1-t^{-1}}\right)^{\#(G)}. 
\end{eqnarray*}
\par
First we compute the $a$-invariant of $A$. 
\begin{description}
\item[Claim 1] ~ 
\begin{enumerate}
 \item (a),(b)$\Longrightarrow [H_{\frm}^d(A)]_j = 0$ for all $j = -1,\ldots,q-d$. 
 \item $\Longrightarrow [H_{\frm}^d(A)]_j = 0$ for all $j = -1,\ldots,q-d,q-d-1$.  
\end{enumerate}
\end{description}
\par \vspace{1mm}
Suppose $(a)$, $(b)$ in $(1)$.  
We may assume that $q \le d-1$. 
Then since $a(k[\Gamma_i]) \le q-d \le -1$, 
we have $[H_{\frm_i}^{d-1}(k[\Gamma_i])]_{0} 
= [H_{\frm_i}^{d-1}(k[\Gamma_i])]_{q-d+1} = 0$ 
for all $i=1,\ldots,n$, where $\frm_i$ denotes the homogeneous 
maximal ideal of $k[\Gamma_i]$. 
\par 
Now let $F$ be a face of $\Delta$ with $1 \le \#(F) \le d-q$. 
As $F$ contains a vertex of $\Delta$ (say $x_i$), 
if we put $G = F \setminus \{x_i\}$, then 
$G \in \Gamma_i$ and $\link_{\Gamma_i} G = \link_{\Delta} F$.  
If $G \ne \emptyset$, then $1 \le \#(G) = \#(F) -1 \le d-q-1$. 
Then 
\[
 \widetilde{H}_{d-\#(F)-1}(\link_{\Delta} F; k) = 
 \widetilde{H}_{d-\#(G)-2}(\link_{\Gamma_i} G; k) = 0  
\]
because $[H_{\frm_i}^{d-1}(k[\Gamma_i])]_{q-d+1} = 0$. 
If $G = \emptyset$, then $F = \{x_i\}$. 
Thus
\[
  \widetilde{H}_{d-\#(F)-1}(\link_{\Delta} F; k) = 
  \widetilde{H}_{d-2}(\Gamma_i; k) = 
  [H_{\frm_i}^{d-1}(k[\Gamma_i])]_{0} =0. 
\]
Hence $\widetilde{H}_{d-\#(F)-1}(\link_{\Delta} F; k) =0$ 
for all $F \in \Delta$ with $1 \le \#(F) \le d-q$. 
This yileds that $[H_{\frm}^d(k[\Delta])]_j = 0$ by Hochster's formula.  
\par 
Suppose $(2)$.  
Then since $a(k[\Gamma_i]) = q-d-1 \le -1$, one can also prove 
the claim in this case by the similar argument as above. // 

\vspace{1mm}
\begin{description}
\item[Claim 2] $[H_{\frm}^d(A)]_0 \cong \widetilde{H}_{d-1}(\Delta;k) =0$. 
\end{description}
\par 
First suppose $(a)$, $(b)$ in $(1)$. 
If $q = d$, then the assertion is clear by the assumption. 
So we may assume that $q \le d-1$. 
Let $K_A$ be the graded canonical module of $A$, that is, 
$[K_A]_j = \Hom_k([H_{\frm}^d(A)]_{-j},k)$.  
Then $[K_A]_1 = 0$ by Claim 1.  
Thus 
\[
 [K_A]_0 \subseteq \bigcap_{i=1}^n (0) :_{K_A} x_i 
= \Hom_A(A/\frm,K_A) = 0, 
\]
where the last vanishing follows from $\depth K_A >0$. 
Thus $[H_{\frm}^d(A)]_{0} = 0$, as required. 
In the case of (2), one can also prove the claim by the same argument as above.  //
\par \vspace{1mm}
By virtue of the above two claims, we get 
\begin{equation} \label{eq:a-inv}
 (1)\;\; a(A) \le q-d-1; \qquad (2)\;\; a(A) \le q-d-2. 
\end{equation}
In particular, in the case of $(2)$, 
we have $\reg A \le a(A) + d+1 \le q-1$ by Hoa--Miyazaki theorem. 
On the other hand, $\reg A \ge \indeg A -1 = q-1$. 
Therefore $A$ has $q$-linear resolution and $a(A) = q-d-2$. 
\par \vspace{1mm}
In order to prove that $A$ has $q$-linear resolution in the case of $(1)$, 
it suffices to show that 
$H_{\frm}^p(A) =0$ for all $p=q,q+1,\ldots,d-1$; 
see Theorem \ref{MinRed}. 
Note that $H_{\frm}^q(A) \cong \widetilde{H}_{q-1}(\Delta;k) =0$ 
by the assumption. 
Also, we have that $\reg A \le q$ by Eq.(\ref{eq:a-inv}). 
Hence $H_{\frm}^p(A) = [H_{\frm}^p(A)]_0 =0$ 
for all $p=q+1,\ldots,d-1$, 
as required.   
The converse follows from Hibi's criterion. 
\end{proof}

\vspace{3mm}
\begin{thm}[\textbf{Characterization of Buchsbaum complex with
\lq\lq minimal multiplicity of type $q$''}]  
\label{Char}
Let $A:=k[\Delta]$ be a $d$-dimensional Buchsbaum Stanley--Reisner ring 
such that $\codim A = c$ and $\indeg A = q$. 
Let $\Delta^{*}$ denote the Alexander dual of $\Delta$ 
and put
\[
 h_{c,d,q} = \frac{(c+q-2)\cdots (c+1)c}{d(d-1)\cdots (d-q+2)}. 
\]
Then the following conditions are equivalent$:$ 
\begin{enumerate}
 \item $A$ has minimal multiplicity of type $q$, that is,  
\[
  e(A) = \frac{c+d}{d} \bbinom{c+q-2}{q-2}.
\] 
\item $A$ has $q$-linear resolution and 
$\dim_k \widetilde{H}_{q-2}(\Delta;k) = h_{c,d,q}$.
 \item The $h$-vector of $A$ is 
\[ 
  \left(1,\,c,\,\cdots, \sbinom{c+q-2}{q-1},\,- \sbinom{d}{q}h,
  \,\sbinom{d}{q+1}h,\cdots,(-1)^{d-q+1} \sbinom{d}{d}h\right),
\]
where $h = h_{c,d,q}$. 
 \item $k[\link_{\Delta} \{x_i\}]$ has $(q-1)$-linear resolution
  for all $i$. 
 \item $a(k[\link_{\Delta} \{x_i\}]) = q-d-1$ for all $i$. 
 \item $a(A) = q-d-2$. 
 \item $k[\Delta^{*}]$ is Cohen--Macaulay with 
 pure and almost linear resolution and with $a(k[\Delta^{*}]) = 0$, that is, 
 the graded minimal free resolution of $k[\Delta^{*}]$ over 
 $S = k[x_1,\ldots,x_n]$ $(n = c+d)$ can be written as follows$:$
\[
 0 \to S(-(c+d))^{\beta_{q}^{*}} \to  
S(-(c+q-2))^{\beta_{q-1}^{*}} \to \cdots \to S(-c)^{\beta_{1}^{*}} 
 \to S \to k[\Delta^{*}] \to 0.    
\]
\end{enumerate}
When this is the case, $\beta^{*}_q = h_{c,d,q}$.  
\end{thm}

\begin{proof}[\quad Proof]
It suffices to show the following implications$:$
$(1)\Rightarrow (2)\Rightarrow (3)\Rightarrow (1)$, 
\noindent 
$(1)\Rightarrow (5)$, $(4) \Leftrightarrow (5)\Leftrightarrow (6)$,  
and $(5),(6)\Rightarrow (7)\Rightarrow (2)$. 
\par \vspace{1mm}
We first show that $(1)\Rightarrow (2)\Rightarrow (3)\Rightarrow (1)$. 
If we suppose (1), then $A$ has $q$-linear resolution
by Theorem \ref{Main}. 
Putting $h=\dim_k \widetilde{H}_{q-2}(\Delta;k)$, we obtain that 
by Proposition \ref{LinMul},  
\[
 \frac{c+d}{d} \bbinom{c+q-2}{q-2} = e(A) 
= \bbinom{c+q-2}{q-2} - h \bbinom{d-1}{q-1}. 
\]
This implies that $h = h_{c,d,q}$. 
In particular, we get (2). 
Also, $(2)\Rightarrow (3)$ follows from Theorem \ref{Key}. 
If we suppose (3), then we have 
\begin{align*}
 e(A) 
 &= \sum_{i=0}^{d} h_i =
    \sum_{i=0}^{q-1} \bbinom{c+i-1}{i} 
   + (-1)^{q+1} \sum_{i=q}^{d} (-1)^{i} \bbinom{d}{i} h_{c,d,q} \\
 & = \bbinom{c+q-1}{q} - \bbinom{d-1}{q-1} h_{c,d,q} 
 = \frac{c+d}{d} \bbinom{c+q-2}{q}. 
\end{align*}
Hence we get (1). 
\par 
Next we show that $(1)\Rightarrow (5)$, 
$(4)\Leftrightarrow (5)\Leftrightarrow (6)$. 
$(1)\Rightarrow (5)$ follows from the proof of Theorem \ref{Main}. 
Since $k[\Gamma_i]$ is Cohen--Macaulay, $(4)$ and $(5)$ are equivalent. 
$(5) \Rightarrow (6)$ follows from Theorem \ref{AinLink},  
and one can prove the converse similarly. 
\par
To complete the proof, we must show that 
$(5),(6)\Rightarrow (7)\Rightarrow (2)$. 
Suppose that $(5)$ and $(6)$. 
By Theorem \ref{AinLink}, $A$ has $q$-linear resolution. 
Thus by Theorem \ref{LinAlex}, the Alexander dual $k[\Delta^{*}]$ is 
a Cohen--Macaulay homogeneous $k$-algebra with almost linear resolution. 
Namely, the graded minimal free resolution of $k[\Delta^{*}]$ over 
a polynomial ring $S$ can be written as follows$:$
\[
 0 \to F_q \to F_{q-1} 
 = S(-(c+q-2))^{\beta_{q-1}^{*}} \to \cdots \to F_1 = S(-c)^{\beta_{1}^{*}} 
 \to S \to A^{*} \to 0,   
\]
where $F_{q} = S(-(c+d))^{\beta^{*}} \oplus S(-(c+q-1))^{\beta^{*'}}$ and  
$\beta^{*} = \dim_k H_{\frm}^{q-1}(A)$,  
$\beta^{*'} = \dim_k H_{\frm}^d(A)_{q-d-1}$.  
Since $a(A) = q-d-2$ by $(6)$, we have $\beta^{*'} =0$, that is, 
the above resolution is pure and thus $a(k[\Delta^{*}]) = 0$. Hence we get (7). 
Conversely, suppose (7). 
By Theorem \ref{LinAlex}, $A=k[\Delta]$  has $q$-linear resolution. 
On the other hand, since $k[\Delta^{*}]$ is a Cohen--Macaulay homogeneous 
$k$-algebra with pure resolution 
of type $(c_1,\ldots,c_q) = (c,c+1,\ldots,c+q-2,c+d)$, we have 
\[ 
 \beta^{*}_q = (-1)^{q+1} \prod_{j=1}^{q-1} \frac{c_j}{c_j-c_q} = h_{c,d,q}.  
\]
by Herzog--K\"uhl's formula (see \cite{HeKu} or \cite[Theorem 4.1.15]{BrHe}).  
Combining with $h = \beta^{*}_q$, we obtain that $h=h_{c,d,q}$, as required. 
\end{proof}

\par\vspace{1mm}
In the following, we give some examples of Buchsbaum Stanley--Reisner rings 
with minimal multiplicity of type $3$. 
We say that a simplicial complex $\Delta$ is spanned by a set $S$ 
if $S$ is the set of facets of $\Delta$. 

\begin{exam}[{\rm Hibi \cite{Hi1} }] \label{Hibi-Ex}
Let $d \ge 2$ be an integer, and let $k$ be a field. 
Put $n = 2d-1$ and $V = \{1,2,\ldots,n\}$. 
Let $\Delta$ be the simplicial complex 
which is spanned by 
$S = \left\{\{\overline{i},\overline{i+1},\ldots, \overline{i+d-1}\}
\,|\, i = 1,2,\ldots, 2d-1\right\}$, 
where $\overline{p}$ stands for $q \in V$ with $p \equiv q \pmod{2d-1}$. 
\par
Then $k[\Delta]$ is a $d$-dimensional Buchsbaum Stanley--Reisner ring with 
minimal multiplicity of type $3$. 
\end{exam}

\begin{proof}[\quad Proof]
Put $A = k[\Delta]$ and $c = n-d = d-1$. 
First note that $\{i,j\} \in \Delta$ for any $i$, $j$ with $1 \le i < j \le n$. 
On the other hand, $\{1,2,d\} \notin \Delta$. 
Thus $A$ is a $d$-dimensional equidimensional Stanley-Reisner ring 
with $\indeg A = 3$. 
Also, we have 
\[
\link_{\Delta} \{1\}  = \{2,\ldots,d\} \cup \{2d-1,2,\ldots,d-1\} 
\cup \cdots \cup \{d+1,d+2,\ldots,2d-1\}. 
\]
Hence $\link_{\Delta} \{1\}$ is a Cohen--Macaulay complex 
which is a $(d-2)$-tree (e.g. cf. \cite{Te3}) and 
thus it has $2$-linear resolution (\cite{Fr2}). 
Similarly, $\link_{\Delta} \{i\}$ is Cohen--Macaulay with 
$2$-linear resolution. 
By Theorem \ref{Char}(4), 
$k[\Delta]$ is a Buchsbaum ring with minimal multiplicity of type $3$. 
\end{proof}

\begin{exam}[\textbf{The Alexander dual of a cyclic polytope}] 
\label{Cyclic-MM}
Let $q$, $d$ be integers with $2 \le q \le d$. 
Put $n = 2d-q+2$ and $f = 2(d-q+1)$. 
Let $C(n,f)$ be a cyclic polytope with $n$ vertices.  
Also, let $\Delta$ be the Alexander dual of 
the boundary complex of $C(n,f)$. 
Then $k[\Delta]$ is a $d$-dimensional 
Buchsbaum Stanley--Reisner ring with minimal multiplicity of type $q$. 
\par
In particular, Conjecture $\ref{Conj-1}$ is true for $h_{c,d,q} =1$. 
\end{exam}

\begin{proof}[\quad Proof]
It follows from Example \ref{CyclicPoly} and Theorem \ref{Char}(2). 
\end{proof}

\begin{exam}[{\rm \cite[Theorem 3.3]{Te1}}] \label{Terai-Ex}
Let $n$ be an integer such that $n >3$, and 
suppose that $2n+1$ is a prime number.  
Let $\Delta$ be the simplicial complex on $V = \{1,2,\ldots,n\}$ 
which is spanned by 
\begin{align*}
 S = & \left\{\{a,b,a+b\}\,|\,1 \le a < b,\; a+b \le n\right\} \\
   & \cup \left\{\{a,b,c\}\,|\, 1 \le a < b < c \le n,\; 
     a+b+c = 2n+1 \right\}.
\end{align*}
Then $A=k[\Delta]$ is a $3$-dimensional 
Buchsbaum Stanley--Reisner ring such that 
$e(A) = \frac{n(n-2)}{3}$. 
In particular, $A$ has minimal multiplicity of type $3$, and thus 
has $3$-linear resolution.  
\end{exam}

\begin{proof}[\quad Proof] 
By \cite{Te1}, $A$ is Buchsbaum with $e(A) = \frac{n(n-2)}{3}$. 
If one put $c = n-3$ and $d=q=3$, then   
$e(A) = \frac{c+d}{d} \bbinom{c+q-2}{q-2}$.
Thus $A$ has minimal multiplicity of type $3$ by definition.
\end{proof}

\par 
In the above example, the assumption \lq\lq $2n+1$ is prime'' is superfluous 
in some sense. 
Indeed, if $n \equiv 0$ or $\equiv 2$ $\pmod{3}$, then 
there exists a $3$-dimensional Buchsbaum Stanley--Reisner ring $k[\Delta]$
with minimal multiplicity which satisfies $\emb(k[\Delta]) =n$; 
see Example \ref{Hanano-Ex} for details. 
\par 
In general, $h_{c,d,q}$ is \textit{not} an integer. 
For example, if $d =q=3$, then since $h_{c,d,q} = \frac{c(c+1)}{6}$. 
However, Hanano constructed Buchsbaum complexes which satisfies 
$h= \dim_k H_{\frm}^2(A) =  \lfloor h_{c,3,3} \rfloor$. 
Motivated by his work, we introduce the following notion.

\begin{defn}[\textbf{Maximal homology}] 
\label{MaxHom}
Let $A=k[\Delta]$ be a $d$-dimensional Buchsbaum Stanley--Reisner ring
with $q$-linear resolution.  
Put $c=\codim A$ and $h = \dim_k H_{\frm}^{q-1}(A)$. 
We say that $A$ has {\it maximal homology\/}
if $h = \lfloor h_{c,d,q}\rfloor$ holds, where 
\[
  h_{c,d,q} =  \frac{(c+q-2)\cdots (c+1)c}{d(d-1)\cdots (d-q+2)}.
\]
\end{defn}

\par 
According to Theorem \ref{Char}, $A$ has minimal multiplicity of type $q$ 
if and only if $A$ has $q$-linear resolution of maximal homology and 
$h_{c,d,q}$ is an integer. 

\begin{exam}[{\rm Hanano \cite{Ha}}] \label{Hanano-Ex}
Let $n$ be an integer with $n \ge 5$. 
Let $\Delta$ be the simplicial complex on $V$ which is spanned by 
the following set $S$. 
Then $k[\Delta]$ is a $3$-dimensional Buchsbaum Stanley--Reisner ring 
with $3$-linear resolution of maximal homology. 
Furthermore, $k[\Delta]$ has minimal multiplicity of type $3$ 
if and only if $n \equiv 0$ or $\equiv 2$ $\pmod{3}$. 
\begin{enumerate}
\item The case of $n \not \equiv 1 \pmod{3}$. 
Put $V:= \{0,1,\ldots,n-1\}$ and  
\begin{eqnarray*}
 S &\! =\! & \left\{\{\overline{i},\overline{i+k},\overline{i+2k}\}\,|\, 
0 \le i \le k-1\right\} \\
   & \quad & \cup \left\{\{\overline{i},\overline{i+k},\overline{i+j}\}\,|\, 
0 \le i \le 3k-1,\,k+1 \le j \le 2k-1\right\}, 
\end{eqnarray*}
when $n = 3k$ and 
\[
 S  = \left\{\{\overline{i},\overline{i+1},\overline{i+3j+2}\}\,|\, 0 \le i \le 3k+1,\,
   0 \le j \le k-1\right\}, 
\]
where $n =  3k+2$. 
Here $\overline{p}$ stands for $q \in V$ with $p\equiv q \pmod{n}$. 
\item The case of $n \equiv 1 \pmod{3}$. 
Put $V := \{\infty,0,1,\ldots,n-2\}$ and 
\begin{eqnarray*}
 S & = &\left\{\{\infty, \overline{i},\overline{i+1}\}\,|\, 
0 \le i \le 3k-1\right\} \\
   & \quad  & \cup \left\{\{\overline{i},\overline{i+1},\overline{i+3j}\} 
\,|\, 0 \le i \le 3k-1,\,1 \le j \le k-1\right\},
\end{eqnarray*}
when $n = 3k+1$. 
Here $\overline{p}$ stands for $q \in V$ with $p\equiv q \pmod{n-1}$. 
\end{enumerate}
\end{exam}

\begin{proof}
(1) By \cite{Ha}, $k[\Delta]$ is a Buchsbaum ring with 
$\indeg k[\Delta] = 3$ and  
\begin{eqnarray*}
(1) \;h(\Delta) & = &  
\left(1,n-3,\frac{(n-2)(n-3)}{2},-\frac{(n-2)(n-3)}{6}\right); \\
(2) \;h(\Delta)  & =& 
\left(1,n-3,\frac{(n-2)(n-3)}{2},-\frac{(n-1)(n-4)}{6}\right).
\end{eqnarray*}
In particular, 
\[
 e(k[\Delta]) = 
\begin{cases}
 \dfrac{n(n-2)}{3} & \text{in (1)}; \\[2mm]
 \dfrac{(n-1)^2}{3} & \text{in (2)}. 
\end{cases}
\]
Therefore the assertion follows from the lemma below. 
\end{proof}

\begin{lemma} \label{3dimMaxHom}
Let $A = k[\Delta]$ be a $3$-dimensional Buchsbaum Stanley--Reisner ring 
of $\Delta$ on $V = \{x_1,\ldots,x_n\}$ 
with $\indeg k[\Delta] = 3$. 
Then 
\begin{enumerate}
 \item $k[\Delta]$ has minimal multiplicity of type $3$ if and only if 
$e(k[\Delta]) = \frac{n(n-2)}{3}$. 
 \item $k[\Delta]$ has maximal homology which is not minimal multiplicity 
of type $3$ if and only if $e(k[\Delta]) = \frac{(n-1)^2}{3}$.
\end{enumerate}
\end{lemma}

\begin{proof}[\quad Proof]
Since $(1)$ just follows from the definition, it suffices to consider 
(2) only. 
\par
Suppose that $e(k[\Delta]) = \frac{(n-1)^2}{3}$. 
Then if we prove that $A = k[\Delta]$ has $3$-linear resolution, 
then by Proposition \ref{LinMul} we have 
\[
\dim_k H_{\frm}^2(A) 
 =  \bbinom{n-1}{2} - e(A) 
=  \frac{(n-1)(n-4)}{6} = \lfloor h_{c,3,3}\rfloor. 
\]
Hence $A$ has maximal homology. 
Thus it is enough to show that $A$ has $3$-linear resolution. 
\par
Put $\Gamma_i = \link_{\Delta} \{x_i\}$ for all $i=1,\ldots,n$. 
By the similar argument as in the proof of Theorem \ref{Main}, 
we have 
\[
 (n-1)^2 = 3 \cdot e(A)= \sum_{i=1}^{n} e(k[\Gamma_i]) \ge n(n-2). 
\]
Note that $e(k[\Gamma_i]) \ge n-2$ and the equality holds if and only if 
$k[\Gamma_i]$ has $2$-linear resolution. 
So we may assume that 
\[
 e(k[\Gamma_1]) = n-1, \qquad e(k[\Gamma_i]) = n-2 \;\;\text{for all $i\ge 2$}. 
\]
Then $a(k[\Gamma_i]) = -1$ and thus $x_i [K_A]_0 =0$ for all $i=2,\ldots,n$. 
Thus 
\[
 [K_A]_0 \subseteq \bigcap_{i=2}^n (0):_{K_A} x_i 
= \Hom_A(A/(x_2,\ldots,x_n)A,K_A) = 0,
\]
where the last vanishing follows from $(x_2,\ldots,x_n)A \notin \Ass_A K_A 
= \Assh(A)$. 
Hence $\widetilde{H}_2(\Delta;k) \cong  [H_{\frm}^3(A)]_0 = 0$. 
Thus $A$ has $3$-linear resolution, as required. 
\end{proof}

\begin{remark} \label{Uniq}
There are two different Buchsbaum complexes with minimal 
multiplicity which has the same numerical data. 
For instance, if $n \ge 8$, $n \equiv 2 \pmod{3}$ and $2n+1$ is prime, 
then $\Delta$ in Example $\ref{Terai-Ex}$ is different 
from $\Delta$ in Example $\ref{Hanano-Ex}$ since 
they have distinct links.  
\end{remark}

\par 
The next example gives Buchsbaum Stanley--Reisner rings 
with minimal multiplicity of higher type.  

\begin{exam}[{\rm Bruns--Hibi \cite[Proposition 3.2]{BrHi}}]
Let $n \ge 6$ be an even integer, and let $k$ be a field. 
Let $\Gamma$ be the simplicial complex whose facets are 
\[
 \{\overline{i},\overline{i+1},\overline{i+2}\},\;
 \{\overline{i},\overline{i+1},\overline{i+4}\},
\ldots,
 \{\overline{i},\overline{i+1},\overline{i+(n-2)}\},\quad
i=1,\ldots,n,
\]
where $\overline{p}$ stands for $q \in [n]$ with $p \equiv q \pmod{n}$. 
Also, let $\Delta$ be the Alexander dual of $\Gamma$. 
Then $A = k[\Delta]$ is a $(n-3)$-dimensional Buchsbaum Stanley--Reisner 
ring with minimal multiplicity of type $(n-3)$. 
Also, $\codim A =3$ and $\emb(A) = n$. 
\end{exam}

\begin{proof}[\quad Proof]
$k[\Gamma]$ is a $3$-dimensional 
Cohen-Macaulay Stanley--Reisner ring with pure and almost 
$3$-linear resolution by \cite[Proposition 3.2]{BrHi}. 
Thus the required assertion follows from Theorem \ref{Char}(7). 
\end{proof}

\vspace{2mm}
\section{Cohen--Macaulay cover of a Buchsbaum complex}

\par
In this section, 
we introduce the notion of \textit{Cohen--Macaulay cover}, 
and prove that such a cover always 
exists for any Buchsbaum simplicial complex $\Delta$ 
with $d$-linear resolution ($d = \dim k[\Delta]$); 
see Theorem \ref{CoverExis}.  

\par 
Also, we prove that if a complex $\Delta$ is 
between two Buchsbaum complexes with $d$-linear resolution 
on the same vertex set $V$, it is also $d$-linear Buchsbaum on $V$. 
Using these facts, we can reduce our problem to 
construct \lq\lq Buchsbaum simplicial complexes with maximal homology'' with 
given parameters.  
As an application, in the case of $d=3$, we will prove that 
Conjecture \ref{Conj-1} is true; see Theorem \ref{P-Ans}.  

\par \vspace{1mm}
In the following, let $\Delta$ be a $(d-1)$-dimensional simplicial complex 
on $V = \{x_1,\ldots,x_n\}$ over a field $k$ with $\indeg k[\Delta] = d$. 

\begin{defn}[\textbf{Cohen--Macaulay cover}] \label{CoverDef}
A simplicial complex $\covDelta$ on $V$ is said to be a  
\textit{Cohen--Macaulay $($$d$-linear$)$ cover} of 
$\Delta$ over $k$ if $\covDelta$ 
satisfies the following conditions$:$ 
\begin{enumerate}
 \item $\covDelta$ is a $(d-1)$-dimensional complex which contains 
  $\Delta$ as a subcomplex. 
 \item $k[\covDelta]$ is a Cohen--Macaulay ring with $d$-linear resolution. 
\end{enumerate}
\end{defn}

\begin{thm}[\textbf{Existence of Cohen--Macaulay cover}] \label{CoverExis}
If $k[\Delta]$ is a $d$-dimensional 
Buchsbaum Stanley--Reisner ring  
with $d$-linear resolution, 
then there exists a Cohen--Macaulay 
$($$d$-linear$)$ cover $\covDelta$ of $\Delta$.  
\end{thm}

\begin{proof}[\quad Proof]
We may assume that $k$ is infinite. 
We use the following notation$:$
\begin{enumerate}
 \item[$\bullet$] $S = k[x_1,\ldots,x_n]$, $A = k[\Delta] = S/I_{\Delta}$
 \item[$\bullet$] $h = \dim_k H_{\frm}^{d-1}(A) = 
 \dim_k \widetilde{H}_{d-2}(\Delta;k)$, $c = n-d$. 
 \item[$\bullet$] $\ecm = \bbinom{n-1}{d-1}$, 
  $\mucm = \bbinom{n-1}{d}$.
 \item[$\bullet$] $\F = \{\{x_{i_1},\ldots,x_{i_d}\}\subseteq V\,|\, 
  \{x_{i_1},\ldots,x_{i_d}\} \;\text{is a facet of $\Delta$}\}$. 
 \item[$\bullet$] $\G = \{G \subseteq V \,|\, \#(G) = d,\, G \in \F\}$. 
\end{enumerate}
Then we note that 
\[
 e := e(A) = \#\F = \ecm - h,\qquad 
 \mu := \mu(I_{\Delta}) = \#\G   = \mucm + h 
\]
by Proposition \ref{LinMul}. 
\par
Take a homogeneous ideal $J \subseteq S$ such that $JA$ is a minimal reduction 
of $\frm A = (x_1,\ldots,x_n)A$. 
If necessary, we may assume that 
$J$ is generated by the following elements $f_{c+1},\ldots,f_n:$
\[
 f_i = x_i - \sum_{j=1}^c c_{i,j}x_j \quad 
(c_{i,j} \in k,\; i=c+1,\ldots,n).
\]
\par
First, we show the following claim$:$ 
\begin{description}
 \item[Claim] $I_{\Delta} + J = (x_1,\ldots,x_c)^d +J$ in $S$. 
\end{description}
\par \vspace{1mm}
Actually, since $A$ is Buchsbaum and $JA$ is a parameter ideal of $A$, 
we have 
\[
  l_S(S/I_{\Delta}+J) = l_A(A/JA) = e(JA) +I(A) =e(A)+h = \ecm.
\]
Also, since $A$ has $d$-linear resolution, by Theorem \ref{MinRed}, 
we have 
\[
  (x_1,\ldots,x_c)^d + J \subseteq I_{\Delta}+J. 
\]
On the other hand, since 
\[
   l_S(S/(x_1,\ldots,x_c)^d + J) = \bbinom{c+d-1}{d-1} = \ecm 
   = l_S(S/I_{\Delta}+J), 
\]
we obtain the required equality, and the claim is proved. 
\par 
Now consider any term order $<$ on $S$, and put 
\[
   \{M_1,\ldots,M_{\mu}\} 
  = \left\{ x_{i_1}\cdots x_{i_d}\,\bigg|\, 
\{x_{i_1},\ldots,x_{i_d}\} \in \G\right\},
\]
where $M_1 < \cdots < M_{\mu}$. 
Also, we put 
\[
  \{N_1,\ldots, N_{\mu-h}\} 
  = \left\{x_1^{k_1}\cdots x_{c}^{k_c}
  \,\bigg|\, k_1+\cdots + k_c = d,\, k_i \ge 0\right\}.
\]
Let $\widetilde{M_i}$ be the homogeneous polynomial of degreed $d$ given 
by substituting $x_i = \sum_{j=1}^c c_{ij}x_j$ for $M_i$. 
Then $\widetilde{M_i}$ can be written as a linear combination of 
$N_1, \ldots,N_{\mu-h}$$:$
\[
   \widetilde{M_i} = a_{i1}N_1 + \cdots + a_{i,\mu-h}N_{\mu-h}
\]
for all $i = 1,\ldots, \mu$. 
Since $I_{\Delta}+J$ generates $(x_1,\ldots,x_c)^d$ in 
$S/J \cong k[x_1,\ldots,x_c]$, the coefficient 
matrix $\A = (a_{ij}) \in \Mat(k; \mu \times (\mu-h))$ satisfies 
that $\rank \A = \mu-h$. 
In particular, there exist $(\mu-h)$-distinct row vectors of $\A$ 
which are linearly independent over $k$. 
Let $G_1,\ldots,G_h$ be elements of $\G$ 
corresponding the other rows of $\A$, and 
put $\widetilde{\Delta} =\Delta \cup \{G_1,\ldots,G_h\}$. 
Then $\widetilde{\Delta}$ be a $(d-1)$-dimensional 
simplicial complex on $V$ contains $\Delta$. 
Also, if we put $\widetilde{A} = k[\widetilde{\Delta}]$, 
then 
\[
  l_{\widetilde{A}}(\widetilde{A}/J\widetilde{A}) 
 = \dim_k k[x_1,\ldots,x_c]/(x_1,\ldots,x_c)^d = \ecm
\]
and $e(\widetilde{A}) = e(A) +h = \ecm$. 
Hence $\widetilde{A}$ is Cohen--Macaulay. 
Furthermore, by Lemma \ref{MulIneq}, $\widetilde{A}$ has $d$-linear 
resolution, and $\widetilde{\Delta}$ 
becomes a Cohen--Macaulay cover of $\Delta$. 
\end{proof}

\par
Let $\Delta^{\rm{min}}$ be a Buchsbaum simplicial complex on $V$ 
with $d$-linear resolution and 
$\dim_k H_{\frm}^{d-1}(k[\Delta^{\rm{min}}]) = 
\lfloor h_{c,d,d}\rfloor = \lfloor \frac{(c+d-2)\cdots (c+1)c}{d!}\rfloor $. 

\begin{cor} \label{MinCover}
If such a complex $\Delta^{\rm{min}}$ exists, then  there exists a 
Cohen--Macaulay cover $\widetilde{\Delta}$ of $\Delta^{\rm{min}}$.  
\end{cor}

\par \vspace{1mm}
In the following, we want to prove that a complex between 
two Buchsbaum complexes  with $d$-linear resolution is 
also Buchsbaum with $d$-linear resolution. 
If this is true, it enables us to construct many Buchsbaum complexes 
with $d$-linear resolution as an application of 
Cohen--Macaulay cover. 
However, it is troublesome to prove the above claim by an induction 
because the links of a simplicial complex 
with $d$-linear resolution do  
not necessarily have linear resolution; 
see also Theorem \ref{Char}. 
In order to get rid of this difficulty,  
we introduce the notion of $d$-fullness as follows. 

\begin{defn} \label{D-fulldef}
Let $d \ge 2$ be an integer. 
Let $\Delta$ be a $(d-1)$-dimensional simplicial complex on 
$V$. 
Then $\Delta$ is called {\it $d$-full\/} if 
$\Delta$ is pure and it contains all $(d-2)$-faces of $2^{V}$. 
\end{defn}

\begin{remark} \label{D-full-remark}
Let $k$ be any field. 
$\Delta$ is $d$-full if and only if $k[\Delta]$ is a $d$-dimensional 
equidimensional Stanley--Reisner ring with $\indeg k[\Delta] \ge d$. 
\par 
If $k[\Delta]$ is a Buchsbaum Stanley--Reisner ring 
with $d$-linear resolution, 
then $\Delta$ is $d$-full.  
\end{remark}

\begin{lemma} \label{D-full}
Let $\Delta^{-} \subseteq \Delta \subseteq \Delta^{+}$ be simplicial 
complexes on $V$. 
Then the following statements hold. 
\begin{enumerate}
 \item If both $\Delta^{-}$ and $\Delta^{+}$ are Cohen-Macaulay $d$-full 
complexes, then so is $\Delta$.  
 \item If both $\Delta^{-}$ and $\Delta^{+}$ are Buchsbaum $d$-full 
complexes, then so is $\Delta$.  
\end{enumerate}
\end{lemma}

\begin{proof}[\quad Proof]
Since $\Delta^{-}$, $\Delta^{+}$ are $d$-full, so is $\Delta$. 
Indeed, as $\Delta^{-} \subseteq \Delta \subseteq \Delta^{+}$, we have that 
$d-1 = \dim \Delta^{-} \le \dim \Delta \le \dim \Delta^{+} = d-1$. 
Thus $\dim \Delta = d-1$. 
In particular, since $\Delta$ can be written as 
$\Delta = \Delta^{-} \cup \{F_1,\ldots,F_h\}$ for some facets of $\Delta^{+}$. 
Moreover, $\indeg k[\Delta] \ge \indeg k[\Delta^{-}] \ge d$. 
\par \vspace{1mm}
(1) We want to prove that $\Delta$ is Cohen--Macaulay by an induction 
on $d= \dim k[\Delta] \ge 2$. 
First suppose that $d =2$. 
Then since $\Delta^{-}$ is connected and $\Delta$ is a complex on 
the same vertex set $V$, $\Delta$ is also connected. 
That is, $\Delta$ is Cohen--Macaulay. 
\par 
Next suppose that $d \ge 3$. 
Put $\Gamma_{i} = \link_{\Delta} \{x_i\}$,  
$\Gamma^{+}_{i} = \link_{\Delta^{+}} \{x_i\}$,  
$\Gamma^{-}_{i} = \link_{\Delta^{-}} \{x_i\}$
for all $i =1,\ldots,n$. 
Then $\Gamma^{+}_{i}$ and $\Gamma^{-}_{i}$ are Cohen--Macaulay $(d-1)$-full 
complexes. 
Applying the induction hypothesis to $\Gamma^{-}_{i} \subseteq 
\Gamma_{i} \subseteq \Gamma^{+}_{i}$, we obtain that 
each $\Gamma_{i}$ is Cohen--Macaulay. 
Then $\Delta$ is Buchsbaum because 
$\Delta$ is pure and $\Gamma_{i}$ is Cohen--Macaulay for all $i$. 
In particular, $H_{\frm}^{d-1}(k[\Delta]) \cong 
\widetilde{H}_{d-2}(\Delta;k)$ and 
$H_{\frm}^p(k[\Delta]) \cong \widetilde{H}_{p-1}(\Delta;k) =0$ 
for all $p \le d-2$. 
\par 
Now consider the following diagram$:$
\[
 \begin{CD}
 C_{d-1}(\Delta^{-}) @>\partial^{-}_{d-1}>> C_{d-2}(\Delta^{-}) 
 @> \partial^{-}_{d-2}>> C_{d-3}(\Delta^{-}) \\
 @V \tau VV  @V id VV @V id VV  \\
 C_{d-1}(\Delta) @> \partial_{d-1}>> C_{d-2}(\Delta) 
 @> \partial_{d-2}>> C_{d-3}(\Delta),  
\end{CD}
\]
where $\tau$ is injective and the last two vertical maps are identity maps. 
Since $\partial^{-}_{d-2} = \partial_{d-2}$ and 
$\tau$ is injective, we have that 
\[
 0 = \widetilde{H}_{d-2}(\Delta^{-};k) = 
\frac{\Ker \partial^{-}_{d-2}}{\Image \partial^{-}_{d-1}} 
\longrightarrow
\frac{\Ker \partial_{d-2}}{\Image \partial_{d-1}} = \widetilde{H}_{d-2}(\Delta;k)
\]
is surjective. 
This yields that $\widetilde{H}_{d-2}(\Delta;k) =0$ and thus 
$k[\Delta]$ is Cohen--Macaulay, as required. 
\par \vspace{1mm}
(2) Considering the links of each vertex $x_i$, we have 
\[
 \link_{\Delta^{-}} \{x_i\} \subseteq  \link_{\Delta} \{x_i\} \subseteq 
 \link_{\Delta^{+}} \{x_i\}. 
\]
By the assumption, we have that two links of both sides are 
Cohen--Macaulay $(d-1)$-full complexes on $V \setminus \{x_i\}$. 
By (1), $\link_{\Delta} \{x_i\}$ is also Cohen--Macaulay. 
Hence $\Delta$ is Buchsbaum. 
\end{proof}

\begin{thm} \label{Int-d-linear}
Let $\Delta^{-} \subseteq \Delta \subseteq \Delta^{+}$ be simplicial 
complexes on $V$. 
If both $k[\Delta^{-}]$ and $k[\Delta^{+}]$ are Buchsbaum 
Stanley--Reisner rings with $d$-linear  
resolutions, then so is $k[\Delta]$.  
\end{thm}

\begin{proof}[\quad Proof]
Since $k[\Delta^{+}]$ (resp., $k[\Delta^{-}]$) is a Buchsbaum 
ring with $d$-linear resolution, 
$\Delta^{+}$ (resp., $\Delta^{-}$) is a Buchsbaum $d$-full  complex. 
Thus $\Delta$ is Buchsbaum by the above lemma. 
Hence by Hibi's criterion, 
it is enough to show $\widetilde{H}_{d-1}(\Delta;k) = 0$. 
\par
Now consider the following diagram$:$
\[
\begin{CD}
C_{d-1}(\Delta) @>\partial_{d-1}>> C_{d-2}(\Delta) 
@>\partial_{d-2}>> C_{d-3}(\Delta) \\
@V \tau VV  @V id VV @V id VV  \\
C_{d-1}(\Delta^{+}) @>\partial^{+}_{d-1}>> C_{d-2}(\Delta^{+}) 
@>\partial^{+}_{d-2}>> C_{d-3}(\Delta^{+}),  
\end{CD}
\]
where $\tau$ is injective. 
Then 
\[
 \widetilde{H}_{d-1}(\Delta;k) = \Ker(\partial_{d-1}) \hookrightarrow 
 \Ker(\partial_{d-1}^{+}) = \widetilde{H}_{d-1}(\Delta^{+};k) =0
\]
since $k[\Delta^{+}]$ has $d$-linear resolution. 
Hence $\widetilde{H}_{d-1}(\Delta;k) = 0$, as required. 
\end{proof}

\vspace{1mm}
\begin{exam}[\textbf{Real projective plane} {\rm \cite[(5.2),(5.4)]{Hi2}}] 
\label{RealProj}
Let $\Delta$ be the simplicial complex on the vertex set 
$V=[6]$ whose maximal faces are 
$\{1,2,5\}$, $\{1,2,6\}$, $\{1,3,4\}$, 
$\{1,3,6\}$, $\{1,4,5\}$, $\{2,3,4\}$, 
$\{2,3,5\}$, $\{2,4,6\}$, $\{3,5,6\}$ 
and $\{4,5,6\}$. 
\par 
If $\chara k \ne 2$, then $k[\Delta]$ is a Cohen--Macaulay ring with 
$3$-linear resolution. 
On the other hand, if $\chara k = 2$, then 
$k[\Delta]$ is a non-Cohen--Macaulay Buchsbaum ring with $\indeg k[\Delta]=3$, 
but it does not have linear resolution.  
\par 
Put $\Delta' = \Delta \setminus \{4,5,6\}$ (M\"obius band). 
Then $k[\Delta']$ is a Buchsbaum ring with $3$-linear resolution and 
$h:=\dim_k H_{\frm}^2(k[\Delta'])=1$ in any characteristic. 
If $\chara k \ne 2$, then $\Delta$ is a Cohen--Macaulay cover of $\Delta'$. 
On the other hand, if $\chara k=2$, then $\Delta' \cup \{1,4,6\}$ is 
a Cohen--Macaulay cover of $\Delta'$, but $\Delta$ is not.  
\par 
Also, since $0 < h < 1 = h_{3,3,3}=2$, $k[\Delta']$ does not have 
minimal multiplicity of type $3$. 
However, one can easily see that $\Delta'$ cannot contain any 
Buchsbaum complex $\Delta''$ having minimal multiplicity of type $3$. 
\end{exam}

\par \vspace{1mm}
In the following, as an application of the notion of Cohen--Macaulay 
cover, we prove Conjecture 
\ref{Conj-1} in the case of $d=q=3$. 

\begin{thm} \label{P-Ans}
Let $c$, $h$ be integers with $c \ge 1$. 
Then the following conditions  are equivalent$:$
\begin{enumerate}
\item There exists a $3$-dimensional Buchsbaum Stanley--Reisner
 ring $A=k[\Delta]$ with $3$-linear resolution
  such that $\codim A = c$ and $\dim H_{\frm}^{2}(A)=h$.
\item The following inequality holds$:$ 
\[
  0 \le h \le 
 h_{c,3,3} = \frac{(c+1)c}{6}.
\]
\end{enumerate}
\end{thm}

\begin{proof}[\quad Proof]
It is enough to show that if $0 \le h \le \frac{(c+1)c}{6}$ then 
there exists a $3$-dimensional Buchsbaum Stanley--Reisner ring 
$k[\Delta]$ with $3$-linear resolution and 
$\dim H_{\frm}^{2}(k[\Delta])=h$. 
\par
For any positive integer $c$, we have an example of 
Buchsbaum complex $\Delta^{-} := \Delta^{\rm{min}}$ 
with $3$-linear resolution of maximal homology over $k$; 
see Example \ref{Hanano-Ex} due to Hanano.   
That is, $\dim_k H_{\frm}^2(k[\Delta^{-}]) 
= \lfloor \frac{(c+1)c}{6} \rfloor=:h_0$. 
By Theorem \ref{CoverExis}, 
we can take a Cohen--Macaulay cover 
$\Delta^{+} := \widetilde{\Delta^{-}}$ 
of $\Delta^{-}$. 
Then $e(k[\Delta^{+}]) - e(k[\Delta^{-}]) =h_0$. 
For a given $h$, let $\Delta$ be a simplicial subcomplex 
of $\Delta^{+}$ containing $\Delta^{-}$ with 
$(h_0-h)+e(k[\Delta^{-}])$ facets. 
Then $k[\Delta]$ is a Buchsbaum Stanley--Reisner ring with 
$3$-linear resolution and $\dim_k H_{\frm}^2(k[\Delta]) = h$ 
by Theorem \ref{Int-d-linear}, as required.  
\end{proof}


\end{document}